\documentclass[11pt]{article}
 \usepackage{benstyle}
 
\usepackage{float}
\usepackage[font={small}]{caption}
\usepackage{longtable}

\title{Parameter selection and numerical approximation properties of Fourier extensions from fixed data}
\author{Ben Adcock \\ Department of Mathematics \\ Purdue University \\ USA \and Joseph Ruan \\ Department of Mathematics \\ Purdue University \\ USA}

\begin{document}
\maketitle

\begin{abstract}
Fourier extensions have been shown to be an effective means for the approximation of smooth, nonperiodic functions on bounded intervals given their values on an equispaced, or in general, scattered grid.  Related to this method are two parameters.  These are the extension parameter $T$ (the ratio of the size of the extended domain to the physical domain) and the oversampling ratio $\eta$ (the number of sampling nodes per Fourier mode).  The purpose of this paper is to investigate how the choice of these parameters affects the accuracy and stability of the approximation.  Our main contribution is to document the following interesting phenomenon: namely, if the desired condition number of the algorithm is fixed in advance, then the particular choice of such parameters makes little difference to the algorithm's accuracy.  As a result, one is free to choose $T$ without concern that it is suboptimal.  In particular, one may use the value $T=2$ -- which corresponds to the case where the extended domain is precisely twice the size of the physical domain -- for which there is known to be a fast algorithm for computing the approximation.  In addition, we also determine the resolution power (points-per-wavelength) of the approximation to be equal to $T \eta$, and address the trade-off between resolution power and stability.
\end{abstract}

\section{Introduction}\label{s:introduction}
In many problems, one is faced with the task of recovering a smooth function $f: [-1,1] \rightarrow \bbC$ to high accuracy from its pointwise samples on an equispaced, or in general, scattered grid.  This problem is challenging, unless the grid points have a specific distribution, since it is difficult to simultaneously ensure both rapid convergence and numerical stability.  In particular, for equispaced data a result of Trefethen, Platte \& Kuijlaars  states that no stable method can converge faster than root-exponentially in the number of data points \cite{TrefPlatteIllCond}, and that any method with more rapid convergence must necessarily be unstable.  

Nevertheless, it has been widely reported that so-called \textit{Fourier extensions} (also known as \textit{Fourier continuations}) lead to effective methods in practice for reconstructions from equispaced or scattered data \cite{FEStability,BoydFourCont,BoydRunge,brunoFEP,DHFEP,LyonFESVD}.  This was confirmed recently in \cite{FEStability} wherein it was shown that Fourier extensions (henceforth abbreviated to FEs) circumvent the stability barrier of \cite{TrefPlatteIllCond} in a certain sense.  Specifically, they converge down to a finite, but user-controlled, maximal accuracy.  

In the method of FEs a function $f : [-1,1] \rightarrow \bbC$ is approximated by a Fourier series of degree $N$ defined on a larger domain $[-T,T]$, where $T>1$ is a user-controlled parameter.  Suppose now that the number of equispaced data points is equal to $2M+1$ for some $M \geq N$.  As shown in \cite{FEStability,LyonFESVD}, if $M = \eta N$ for some fixed \textit{oversampling} ratio $\eta \geq 1$, then one may compute an accurate and stable FE approximation of $f$ from this data by a simple least-squares fit (see also \S \ref{s:FE} for details).  Moreover, when the extension parameter $T$ is equal to $2$ -- that is, the extended domain $[-T,T]$ is precisely twice the size of the physical domain $[-1,1]$ -- an algorithm developed by M.\ Lyon allows for the computation of the FE approximation in only $\ord{M (\log M)^2}$ operations  \cite{LyonFast}.  Note that this fast algorithm relies on the particular symmetries of FEs when $T=2$.

From a practical standpoint, it is vitally important to understand how to choose the parameters $T$ and $\eta$.  The purpose of this paper is to address this issue.  In particular, we seek to determine how choices of these parameters affect both the stability and the accuracy of the FE approximation.  Note that, due to the aforementioned fast algorithm, there is a seeming advantage to using the value $T=2$.  However, is choice optimal vis-a-vis the other properties of the algorithm, namely, convergence and stability? Or does another choice (albeit lacking a fast algorithm) give better numerical performance in these respects?  We shall provide answers to these questions.

It is known that different choices of $T$ affect the intrinsic approximation properties of the FE approximation system, i.e.\ the space of trigonometric polynomials of degree $N$ on the extended domain $[-T,T]$ \cite{BADHFEResolution,FEStability,BoydFourCont}.  For example, when $T$ is close to $1$, the FE approximation system possesses better \textit{resolution} power for oscillatory functions \cite{BADHFEResolution}. However, it is perhaps not surprising that such choices also require larger oversampling parameters $\eta$ to maintain the algorithm's stability.  When $M$ is fixed, a larger $\eta$ means a smaller parameter $N$, and therefore the best approximation error in the above subspace, which is determined by the size of $N$, is correspondingly larger.  

From this argument, it is apparent that a balance must be struck between $T$ and $\eta$ so as to preserve accuracy and stability.  The main result we obtain in this paper through numerical experiment is that these two effects precisely counteract each other.  Specifically, if the desired condition number of the algorithm is fixed in advance, then, provided $T$ is not too large, the precise choice of $T$ makes no substantial difference to the accuracy of the Fourier extension algorithm.  Smaller $T$ is exactly offset by the requirement of a larger value of $\eta$ to preserve stability.

From this result, we are able to draw two main conclusions.  First, any attempt to optimize $T$ will only bring limited, and most likely highly function-dependent, success.  Second, since the choice of $T$ makes little difference, one may safely use $T=2$, and the resulting fast algorithm, without worrying that this choice may be suboptimal.   We remark in passing that some previous insight into the effect of the parameters was given by Bruno et al.\ in \cite{brunoFEP}.  However, this was largely carried out for specific functions.

On the face of it, the conclusion we draw may appear surprising, or at the very least, a peculiar phenomenon isolated to the particular choice of equispaced data.  After further numerical experiments, we conclude that this phenomenon is actually quite widespread.  Specifically, we show exactly the same results for both scattered nonequispaced data, as well as Fourier data.  Hence we conclude that unless the data is chosen specifically to favour a particular choice of $T$ (see \S \ref{s:other_data} for an example of such data), the value of $T$ makes little difference to the algorithm.  We note here that, much as in the case of equispaced data, a fast Fourier extension algorithm for scattered data has also been developed in the case $T=2$ \cite{LyonFENonuniform}.

In some applications, including the numerical solution of PDEs, an important question about an approximation algorithm is that of resolution power.  Specifically, how many measurements (e.g.\ equispaced function samples) are required to recover an oscillation of frequency $\omega$.  This topic was first investigated rigorously by Gottlieb \& Orszag \cite{naspec}, who popularized the concept of \textit{points-per-wavelength}.  Through our experiments we establish that this quantity for FE approximations is given by the product of $T$ and $\eta$.  We give theoretical arguments as to why this should be the case, and in the case $T=2$ (and therefore, by the above discussion, all values of $T$) provide numerical results assessing the tradeoff between resolution power and numerical stability.

As the reader will have noticed, our aim in this paper is to investigate FEs through numerical experiment.  The main conclusion we draw is based solely on the result of these experiments.  Although we do present some mathematical insight as to why it should hold, this is a ways short of a proof.  We leave this as a topic for future work.  As we explain in \S \ref{s:conclusion}, this will likely require an intricate analysis of the singular values and singular vectors of a certain matrix related to Slepian's prolate matrix \cite{SlepianV,Varah}, which is beyond the scope of this paper.  Nevertheless, we feel the conclusion we draw, albeit without a proof, is of substantial independent interest for anyone seeking to use FEs in practice.

\section{Fourier extensions}\label{s:FE}
Our concern in this paper is the approximation of functions defined on compact intervals, which without loss of generality we take to be $[-1,1]$.  The method of Fourier extensions (FEs) is based on approximating such functions using a Fourier series defined on an extended interval $[-T,T]$, where $T>1$ is the so-called extension parameter.  In other words, given $N \in \bbN$ we compute an approximation to $f$ from the subspace
\bes{
\cG^{(T)}_N : = \spn \left \{ \phi_n : |n| \leq N \right \},\qquad \phi_n(x) = \E^{\I \frac{n \pi}{T} x }.
}

\subsection{Approximation properties of the subspace $\cG^{(T)}_N$}
We now present several results concerning the intrinsic approximation properties of the subspace $\cG^{(T)}_{N}$.  We use the notation $\rH^{k}(-1,1)$ for the standard Sobolev space on an interval $(-1,1)$, where $k \geq 0$.  We denote the corresponding norm by $\nm{\cdot}_{\rH^k(-1,1)}$.

\thm{[\cite{BADHFEResolution}]
\label{t:error_alg}
Let $T_0 > 1$ and suppose that $f \in \rH^{k}(-1,1)$ for some $k \geq 0$.  Then, for each $N \in \bbN$ and $T \geq T_0$, there exists a $\phi \in \cG^{(T)}_{N}$ such that
\bes{
\| f - \phi \|_{\rL^2(-1,1)} \leq C(k,T_0) \left ( \frac{N \pi}{T} \right )^{-k} \| f \|_{\rH^k(-1,1)},\qquad \| \phi \|_{\rL^2(-T,T)} \leq C(k,T_0) \| f \|_{\rH^k(-1,1)},
}
for some constant $C(k,T_0)$ depending on $k$ and $T_0$ only.  
}

This theorem asserts \textit{algebraic} convergence of the best approximations in $\cG^{(T)}_{N}$ when $f$ has $k$ derivatives, and \textit{superalgebraic} convergence whenever $f$ is smooth.  Note that it also implies the existence of a function $\phi \in \cG^{(T)}_{N}$ which gives such convergence rates, and which cannot grow too large on the extended domain $[-T,T]$.  This will be of significance in \S \ref{ss:CondErrFE}.

Our next result confirms \textit{geometric} convergence of best approximations in $\cG^{(T)}_{N}$ in the case that $f$ is analytic.  To state this result, we first recall the definition of a Bernstein ellipse:
\bes{
\cB(\rho) = \left \{ \tfrac12 \left ( \rho^{-1} \E^{\I \theta} + \rho \E^{-\I \theta} \right ) : \theta \in [-\pi,\pi] \right \} \subseteq \bbC,\quad \rho > 1.
}
As discussed in \cite{FEStability,DHFEP}, Fourier extensions can be viewed as polynomial approximations in the mapped variable $z = m(x)$, where
\be{
\label{mapping}
m(x) = 2 \frac{\cos \frac{\pi}{T} x - \cos \frac{\pi}{T} }{1- \cos \frac{\pi}{T}} - 1.
}
Note that $m$ maps $[0,1]$ to $[-1,1]$ bijectively.  Since the convergence of polynomial approximations of analytic functions is determined by Bernstein ellipses, it makes sense to introduce the new regions
\bes{
\cD(\rho) = m^{-1} (\cB(\rho)),\quad \rho > 1,
}
We now have the the following theorem:

\thm{[\cite{DHFEP}]
\label{t:error_geo}
Let $T>1$ be given and suppose that $f$ is analytic in $\cD(\rho')$ for some $\rho' > 1$ and continuous on its boundary.  Then, for each $N \in \bbN$, there exists a $\phi \in \cG^{(T)}_{N}$ such that
\be{
\label{geo_err}
\| f - \phi \|_{\rL^{\infty}(-1,1)} \leq \frac{c_f(T)}{1-\rho} \rho^{-N},
}
where $c_f(T) > 0$ is proportional to $\max_{z \in D(\rho) } | f(z) |$,
\bes{
\rho = \min \left \{ \rho' , E(T) \right \},
}
and $E(T) = \cot^2 \left ( \frac{\pi}{4 T} \right )$.  Moreover, $\phi$ satisfies
\be{
\label{geo_norm}
\| \phi \|_{\rL^{\infty}(-T,T)} \leq c_f(T) \left ( E(T) / \rho \right )^N.
}
}

This theorem establishes geometric convergence of best approximations in $\cG^{(T)}_{N}$.  However, \R{geo_norm} suggests that in order to obtain such a convergence rate, one may have to allow for exponential growth of the corresponding $\phi$ in the extended domain $[-T,T]$ whenever $\rho < E(T)$.  We shall return to this observation in \S \ref{ss:CondErrFE}.  

We remark also that Theorem \ref{t:error_geo} asserts that the maximal rate of geometric convergence is limited to $E(T)$, even if $f$ is entire.  This is due to the mapping $m^{-1}$ which introduces a square-root type singularity and thereby limits the overall rate of convergence.  See \cite{FEStability,DHFEP} for a discussion.

\subsection{Fourier extensions from equispaced data}\label{ss:FE_equi}
The concern of the majority of this paper is the approximation of a function $f$ from its values
\bes{
f(m/M),\quad m=-M,\ldots,M,
}
on an equispaced grid of $2M+1$ points.  For convenience, let us define the operator
\bes{
S_M : \rL^{\infty}(-1,1) \rightarrow \bbC^{2M+1},\ f \mapsto \frac{1}{\sqrt{M}}(f(m/M))^{M}_{m=-M}.
}
We refer to $S_M$ as the \textit{sampling} operator.  Given the vector $S_M (f)$ of samples of $f$, we construct its FE approximation in the standard way via a least-squares data fit  \cite{FEStability,BoydFourCont,brunoFEP,DHFEP,LyonFESVD}.  Let $N \leq M$ be given.  Then we define the FE approximation as follows:
\be{
\label{FE_LS_fn}
F^{(T)}_{N,M}(f) : = \underset{\phi \in \cG^{(T)}_N}{\operatorname{argmin}} \sum_{|m| \leq M} \left | f(m/M) - \phi(m/M) \right |^2,
}
or more succinctly,
\bes{
F^{(T)}_{N,M}(f) : =\underset{\phi \in \cG^{(T)}_N}{\operatorname{argmin}} \left | S_M(f-\phi) \right |,
}
where $|\cdot|$ denotes the usual Euclidean norm on $\bbC^{2M+1}$.  Note that $F_{N,M}$ is an operator with domain $\rL^\infty(-1,1)$ and range $\cG^{(T)}_N$.  Moreover, if we denote 
\bes{
F^{(T)}_{N,M}(f) = \sum_{|n| \leq N} a_n \E^{\I \frac{n \pi}{T} x},
}
then the vector $\mathbf{a} = (a_n)^{N}_{n=-N}$ of FE coefficients is the solution of the least squares problem
\be{
\label{FE_LS_coeff}
\mathbf{a} = \underset{\mathbf{c} \in \bbC^{2N+1}}{\operatorname{argmin}} \left |A^{(T)} \mathbf{c} - S_M(f) \right |,
}
where $A^{(T)} \in \bbC^{(2M+1) \times (2N+1)}$ has entries
\bes{
(A^{(T)})_{m,n} = \frac{1}{\sqrt{M}} \E^{\I \frac{n m \pi}{T}},\quad |m| \leq M,\ |n| \leq N.
}
Note that the normalization $1/\sqrt{M}$ in both $S_M$ and $A^{(T)}$ means that the entries of the normal matrix $(A^{(T)})^* A^{(T)}$ are Riemann sum approximations to the Gram matrix of the functions $\phi_n(x)$.  This ensures that the singular values of $A^{(T)}$ lie between $0$ and $1$ for large $M$, which, since we typically consider truncated SVDs with a fixed truncation parameter (see later), ensures that there is no linear drift in the error for large $M$.

For convenience, let us now introduce some additional notation.  Let
\bes{
L^{(T)}_{N,M} : \bbC^{2M+1} \rightarrow \bbC^{2N+1},
}
be defined by
\be{
\label{LNM_def}
L^{(T)}_{N,M}(\mathbf{b}) = \underset{\mathbf{c} \in \bbC^{2N+1}}{\operatorname{argmin}} \left | A^{(T)} \mathbf{c} - \mathbf{b} \right |,
}
and let
\bes{
R^{(T)}_{N} : \bbC^{2N+1} \rightarrow \cG^{(T)}_N,\ \mathbf{a} = (a_n)^{N}_{n=-N} \mapsto \sum_{|n| \leq N} a_n \E^{\I \frac{n \pi}{T} x}.
}
Note that $F^{(T)}_{N,M} = R^{(T)}_{N} \circ L^{(T)}_{N,M} \circ S_M$.

As discussed in \cite{FEStability}, the algebraic least-squares problem to be solved in \R{LNM_def} is highly ill-conditioned.  When applied to \R{LNM_def}, different numerical algorithms may consequently give somewhat different results.  For this reason, it is important to specify the solver used.    In the majority of this paper, as has been previously considered in \cite{FEStability,BoydFourCont,LyonFESVD}, we solve \R{LNM_def} by using truncated singular value decompositions (SVDs).  If $U \Sigma V^*$ denotes the SVD of $A^{(T)}$, where $\Sigma$ is the diagonal matrix of singular values $\sigma_1 \geq \sigma_2 \geq \ldots$, then we correspondingly define
\bes{
L^{(T,\epsilon)}_{N,M}(\mathbf{b}) = V \Sigma^{(\epsilon)} U^* \mathbf{b},
}
where $\Sigma^{(\epsilon)}$ is the diagonal matrix with $n^{\rth}$ entry $1/\sigma_n$ if $\sigma_n > \epsilon$ and $0$ otherwise.  Here $\epsilon > 0$ is the truncation parameter, which we take to be $10^{-13}$ unless specified otherwise.  We denote the corresponding FE by $F^{(T,\epsilon)}_{N,M}(f)$.

Having said this, we note that the quantities introduced below for studying equispaced FEs -- namely, the condition number and numerical defect constant -- are not specific to the SVD algorithm.  In particular, one can compute such quantities for each different numerical solver and thereby directly compare the effectiveness of an equispaced FE resulting from an SVD with an equispaced FE computed using \textit{Matlab}'s $\backslash$ or \textit{Mathematica's} \texttt{LeastSquares} commands, for example.  We return to this briefly in \S \ref{ss:other_solvers}.

\rem{
Regardless of the solver used, it is important that the least-squares \R{LNM_def} is regularized when solved numerically.  This is done by the parameter $\epsilon$ with the SVD approach, or automatically when using an blackbox least-squares solver such as \textit{Matlab}'s $\backslash$ or \textit{Mathematica's} \texttt{LeastSquares}.  

As shown in \cite{FEStability}, the `exact' FE mapping $f \mapsto F^{(T)}_{N,M}(f)$, i.e.\ that obtained by solving \R{LNM_def} in infinite precision, is ill-conditioned and suffers from a Runge phenomenon unless the number of equispaced points $M$ scales quadratically with $N$.  Such severe scaling is undesirable, and is due solely to the behaviour of the Fourier series corresponding to singular vectors with small singular values.  Fortunately, when the system \R{LNM_def} is regularized and $F^{(T,\epsilon)}_{N,M}$ is computed, this scaling drops to linear in $N$.  Moreover, the Fourier series of the excluded singular values are precisely those which are small on the domain $[-1,1]$ but large on $[-T,T] \backslash [-1,1]$.  Thus, their exclusion has little effect on the approximation of $f$.  Note that a similar behaviour is also witnessed when different solvers are used for \R{LNM_def}, such as those listed above.
}

\subsection{Condition number and error bounds for equispaced FE approximations}\label{ss:CondErrFE}

We now provide estimates for the accuracy and stability of $F^{(T)}_{N,M}(f)$.  The key point is that these formulae involve constants which can be computed numerically.  This will be discussed in the next section.  First, however, we require the following assumption:

\vspace{1pc} \noindent \textbf{Assumption.}  The operator $L^{(T)}_{N,M}$ defined by solving \R{LNM_def} with a standard numerical solver (e.g.\ truncated SVDs) is approximately a linear operator.

\vspace{1pc} Note that the exact, i.e.\ infinite precision, version of $L^{(T)}_{N,M}$ is of course a linear operator.  Hence it is not unreasonable that its finite precision counterpart acts in the same way.  Observe also that this assumption implies that the overall numerical FE operator $F^{(T)}_{N,M}$ is also linear.

With this in hand, we can now define the condition number in the usual way:

\defn{[Condition number]
The (absolute) condition number of the equisapced FE approximation $F^{(T)}_{N,M}$ is given by
\be{
\label{kappa_def}
\kappa^{(T)}_{N,M} = \max_{\substack{\mathbf{b} \in \bbC^{2M+1} \\ \mathbf{b} \neq 0}} \left \{ \frac{\nmu{R^{(T)}_N \circ L^{(T)}_{N,M}(\mathbf{b})}_{\infty} }{|\mathbf{b}|_{\infty}} \right \}.
}
Here $\nm{g}_{\infty} = \sup_{x \in [-1,1]} |g(x)|$ is the uniform norm on $[-1,1]$ for $g \in \rL^\infty(-1,1)$ and $|\mathbf{b}|_{\infty} = \max_{|m| \leq M} | b_m|$ for $\mathbf{b} = (b_m)_{|m| \leq M} \in \bbC^{2M+1}$.
}

We remark that $\kappa^{(T)}_{N,M}$ is the absolute condition number, as opposed to the more standard relative condition number \cite{TrefethenBau}.  It measures the absolute sensitivity of the FE to perturbations in the samples of $f$, and transpires to be substantially easier to compute in practice.  Note also that the definition implicitly assumes linearity of the mapping $L^{(T)}_{N,M}$.

We now consider the approximation error.  For this we require the following definition:

\defn{[Numerical defect constant]
The numerical defect constant of the equispaced FE is given by
\be{
\label{lambda_def}
\lambda^{(T)}_{N,M} = \max_{\substack{\mathbf{a} \in \bbC^{2N+1} \\ \mathbf{a} \neq 0}} \left \{ \frac{\nm{R^{(T)}_N \left ( \mathbf{a} - L^{(T)}_{N,M} \circ S_M \circ R^{(T)}_N (\mathbf{a}) \right )}_{\infty}}{| \mathbf{a} |_{\infty}} \right \}.
}
}

Before showing the relevance of this constant to error bounds, let us first consider its meaning.  Recall that each vector $\mathbf{a}$ corresponds uniquely to a function $\phi \in \cG^{(T)}_{N}$ given by $\phi = R^{(T)}_{N} \mathbf{a}$.  Thus the numerator in \R{lambda_def} reads
\bes{
\nm{\phi - F^{(T)}_{N,M}(\phi)}_{\infty}.
}
In infinite precision, the FE operator $F^{(T)}_{N,M}$ satisfies $F^{(T)}_{N,M}(\phi) = \phi$ for $\phi \in \cG^{(T)}_{N}$.  In other words, it is a projection.  Hence the numerical defect constant measures how close the numerical, i.e.\ finite precision, FE operator is to possessing this property.

We are now able to provide an error bound for $F^{(T)}_{N,M}(f)$:

\lem{
\label{l:error_numerical}
Let $f \in \rL^\infty(-1,1)$ and suppose that $F^{(T)}_{N,M}(f)$ is given by \R{FE_LS_fn}.  Then
\be{
\label{error_numerical}
\| f - F^{(T)}_{N,M}(f) \|_{\infty} \leq \inf_{\mathbf{a} \in \bbC^{2N+1}} \left \{ \left ( 1 + \kappa^{(T)}_{N,M} \right )\| f - R^{(T)}_N (\mathbf{a}) \|_{\infty} + \lambda^{(T)}_{N,M} | \mathbf{a} |_{\infty} \right \},
}
where $\lambda^{(T)}_{N,M}$ is as in \R{lambda_def}.
}
\prf{
Let $\mathbf{a} \in \bbC^{2N+1}$ be arbitrary and write $\phi = R^{(T)}_{N}(\mathbf{a})\in \cG^{(T)}_N$.  Then, using linearity of $F^{(T)}_{N,M}$ (Assumption 1), we obtain
\bes{
\| f - F^{(T)}_{N,M}(f) \|_{\infty} \leq \| f - \phi \|_{\infty}+ \| F^{(T)}_{N,M}(f-\phi) \|_{\infty} + \| \phi - F^{(T)}_{N,M}(\phi) \|_{\infty} .
}
We consider the latter two terms separately.  For the first, note that
\eas{
\| F^{(T)}_{N,M}(f-\phi) \|_{\infty} = \| R^{(T)}_N \circ L^{(T)}_{N,M} \circ S_M ( f - \phi) \|_{\infty} \leq  \kappa^{(T)}_{N,M} | S_M(f-\phi) |_{\infty} \leq \kappa^{(T)}_{N,M} \| f - \phi \|_{\infty}.
}
This gives the corresponding second term in \R{error_numerical}.  We now consider the other term.  We have
\bes{
\nm{\phi - F^{(T)}_{N,M}(\phi)}_{\infty} =  \nm{R^{(T)}_N \left ( \mathbf{a} - L^{(T)}_{N,M} \circ S_M \circ R^{(T)}_M (\mathbf{a}) \right )}_{\infty}
\leq \lambda^{(T)}_{N,M} | \mathbf{a} |_{\infty},
}
as required.
}

Let us now interpret this error bound.  In \S \ref{ss:computing} we shall observe numerically that
\be{
\label{kappa_lambda_growth}
\kappa^{(T)}_{M/\eta,M} \sim \tilde{\kappa}_{T,\eta} \log M,\qquad \lambda^{(T)}_{M/\eta,M} \sim \tilde{\lambda}_{T,\eta} M,\qquad M \rightarrow \infty,
}
where $\tilde{\kappa}_{T,\eta}$ and $\tilde{\lambda}_{T,\eta}$ are independent of $M$.  Furthermore, the ratio $\mu = \tilde{\lambda}_{T,\eta}/\tilde{\kappa}_{T,\eta}$ is roughly $10^{-13}$ in magnitude, regardless of the choice of $T$ or $\eta$.  Hence, one has the estimate
\be{
\label{error_numerical_2}
\| f - F^{(T)}_{N,M}(f) \|_{\infty} \leq \tilde{\kappa}_{\eta,T} M \inf_{\mathbf{a} \in \bbC^{2N+1}} \left \{ \| f - R^{(T)}_N (\mathbf{a}) \|_{\infty} + \mu | \mathbf{a} |_{\infty} \right \},\qquad N = M/\eta.
}
The key aspect of this bound is that it separates the error into two component.  The first, namely, $\tilde{\kappa}_{\eta,T}$, is determined by the parameters $\eta = M/N$ and $T$, and is independent of the function $f$.  Moreover, as we see next, it can be computed numerically.  The second, i.e.\ the term
\be{
\label{EN_def}
E_{N}(f) : = \inf_{\mathbf{a} \in \bbC^{2N+1}} \left \{ \| f - R^{(T)}_N (\mathbf{a}) \|_{\infty} + \mu | \mathbf{a} |_{\infty} \right \},
}
is crucially independent of $\eta$ and depends only on the intrinsic approximation properties of the subspace $\cG^{(T)}_N$ and the smoothness of $f$.  In particular, combining Lemma \ref{l:error_numerical} with Theorems \ref{t:error_alg} and \ref{t:error_geo}, we immediately obtain the following:
\cor{
\label{c:EN_rate}
Let $E_N(f)$ be given by \R{EN_def}.  If $f \in \rH^k(-1,1)$ then
\be{
\label{alg_conv}
E_N(f) \leq  \min_{0 \leq l \leq k} \left \{ C(l,T_0) \| f \|_{\rH^l(-1,1)} \left ( \left ( \frac{N \pi}{T} \right )^{-l} + \mu \right ) \right \},
}
for each $T \geq T_0$, where $T_0$ and $C(l,T_0)$ are as in Theorem \ref{t:error_alg}.  Moreover if $f$ is analytic in $\cD(\rho)$ and continuous on its boundary, then one has
\be{
\label{exp_conv}
E_N(f) \leq c_f(T) \rho^{-N} \left (1 + \mu E(T)^N \right ) ,
}
where $c_f(T)$ is as in Theorem \ref{t:error_geo}.
}

This corollary explains the behaviour of $E_N(f)$ in both finite and infinite precision.  In infinite precision, where $\mu = 0$, the bound \R{exp_conv} shows geometric decay of $E_N(f)$ for all $N$ at a rate equal to $\rho$.  In finite precision, however, the small, but nonzero constant $\mu$ dramatically alters the convergence.  Geometric decay still occurs for small $N$, when the term $\mu E(T)^N$ in \R{exp_conv} is small, but once $N \geq N_0 = - \log \mu / \log E(T)$ the right-hand side of \R{exp_conv} begins to increase.  For $N \geq N_0$, $E_N(f)$ no longer decays geometrically.  Instead, its decay is described by the bound \R{alg_conv}.  Specifically, algebraic decay in $N$ occurs down to a maximal achievable accuracy on the order of $\mu$.  Note that the constant term $C(l,T) \| f \|_{\rH^l(-1,1)}$ usually grows with $l$, thus as $E_N(f)$ approaches $\mu$ the effective rate of decay usually lessens.  

\begin{figure}
\begin{center}
$\begin{array}{ccc}
\includegraphics[width=6.00cm]{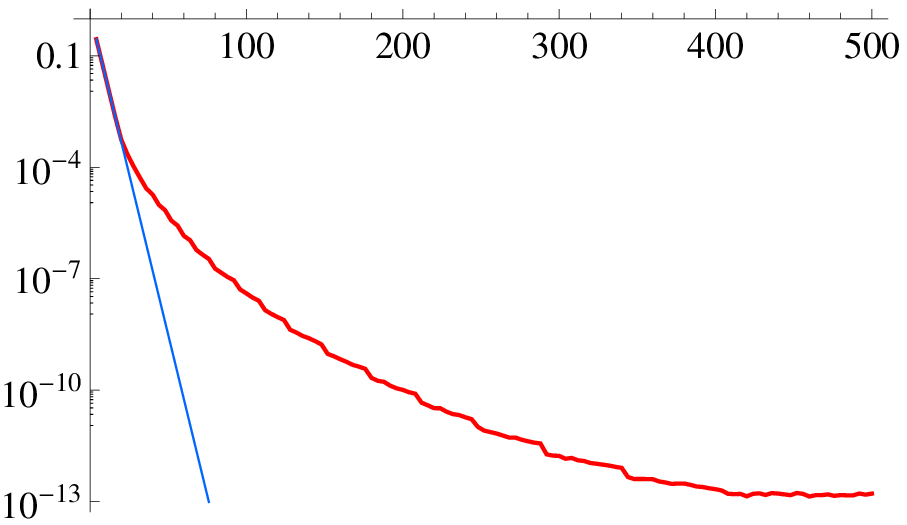}  &  \hspace{2pc} &\includegraphics[width=6.00cm]{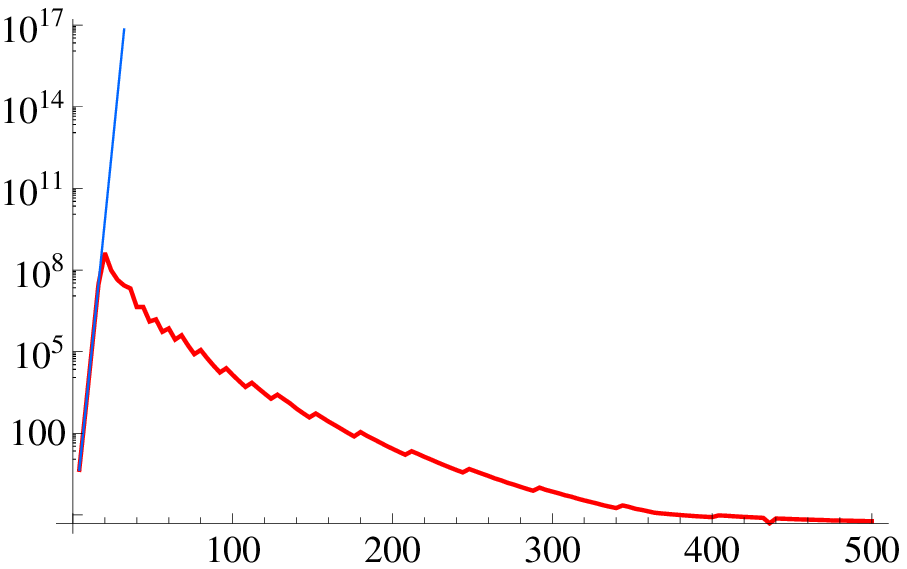} 
\end{array}$
\caption{Best approximation error (left) and coefficient norm (right) for the function $f(x) = \frac{1}{40-39 x}$ in finite precision (thicker line) and infinite precision (thinner line).  The latter was computed in \textit{Mathematica} using additional precision.} \label{f:EN_Plot}
\end{center}
\end{figure}

This behaviour is illustrated in Figure \ref{f:EN_Plot}.  As we see, geometric convergence in infinite precision requires geometric growth of the coefficient vector $\mathbf{a}$.  Conversely, in finite precision, such convergence is sacrificed for algebraic convergence whilst maintaining a bounded coefficient norm.  We refer to \cite{FEStability} for a more detailed discussion.

\subsection{Computing the condition number and numerical defect constant}\label{ss:computing}
Whilst Corollary \ref{c:EN_rate} explains the decay of $E_N(f)$, in order to understand the error of the equispaced FE $F^{(T)}_{N,M}(f)$ we need to determine the magnitudes of $\kappa^{(T)}_{N,M}$ and $\lambda^{(T)}_{N,M}$.  The former also determines the stability of $F^{(T)}_{N,M}$.  As we discuss further in \S \ref{s:conclusion}, it is as of yet unknown how to do this analytically, hence we now resort to numerical investigations.  For this we need a means of computing $\kappa^{(T)}_{N,M}$ and $\lambda^{(T)}_{N,M}$.  This follows from the next two lemmas:

\lem{
\label{l:kappa_equals}
Let $\mathbf{e}_m$, $|m| \leq M$, be the canonical basis for $\bbC^{2M+1}$.  Then
\bes{
\kappa^{(T)}_{N,M} = \sup_{x \in [-1,1]} \sum_{|m| \leq M} \left | R^{(T)}_N \circ L^{(T)}_{N,M}(\mathbf{e}_m)(x) \right |.
}
}
\prf{
If $\mathbf{b} = (b_m)_{|m| \leq M} \in \bbC^{2M+1}$, we may write $\mathbf{b} = \sum_{|m| \leq M} b_m \mathbf{e}_m$.  By linearity
\bes{
R^{(T)}_{N} \circ L^{(T)}_{N,M}(\mathbf{b}) = \sum_{|m| \leq M} b_m R^{(T)}_{N} \circ L^{(T)}_{N,M}(\mathbf{e}_m),
}
and therefore 
\bes{
\kappa^{(T)}_{N,M} \leq \sup_{x \in [-1,1]} \sum_{|m| \leq M} | R^{(T)}_{N} \circ L^{(T)}_{N,M}(\mathbf{e}_m) | .
}
Conversely,
\bes{
\kappa^{(T)}_{N,M} = \sup_{x \in [-1,1]} \max_{\substack{\mathbf{b} \in \bbC^{2M+1} \\ \mathbf{b} \neq 0}} \frac{\left | \sum_{|m| \leq M} b_m R^{(T)}_{N} \circ L^{(T)}_{N,M}(\mathbf{e}_m)(x) \right |}{| \mathbf{b} |_{\infty} }.
}
We now set $b_m$ equal to the complex sign of $R^{(T)}_N \circ L^{(T)}_{N,M}(\mathbf{e}_m)(x)$ to deduce the lower bound.
}
This lemma allows for approximate computation of $\kappa^{(T)}_{N,M}$.  Let $K \in \bbN$ be given and define
\bes{
x_k = \frac{T(k-1)}{K}-1,\quad k=1,\ldots,K_T,
}
where
\bes{
K_T = \left \lfloor \frac{2 K}{T} +1 \right \rfloor .
}
Note that $\{ x_k \}^{K_T}_{k=1}$ is a set of $K_T$ equispaced nodes in $[-1,1]$.  Therefore
\bes{
\kappa^{(T)}_{N,M} = \lim_{K \rightarrow \infty} \kappa^{(T)}_{N,M,K},
}
where
\bes{
\kappa^{(T)}_{N,M,K} = \max_{k=1,\ldots,K_T} \sum_{|m| \leq M} \left | R^{(T)}_N \circ L^{(T)}_{N,M}(\mathbf{e}_m)(x_k) \right |,
}
is a computable quantity.  We remark also $\kappa^{(T)}_{N,M,K}$ can be computed efficiently using Fast Fourier Transforms (FFTs), since the functions $R^{(T)}_N \circ L^{(T)}_{N,M}(\mathbf{e}_m)(x)$ are Fourier series and $\{ x_k \}^{K_T}_{k=1}$ are appropriately constructed equispaced nodes.  Throughout this paper we shall consistently use the value $K=2^{15}$ in our numerical experiments.

We use a similar approach in order to compute the numerical defect constant.  Analogously to Lemma \ref{l:kappa_equals}, we have the following:
\lem{
Let $\mathbf{e}_n$, $|n| \leq N$, be the canonical basis for $\bbC^{2N+1}$.  Then
\bes{
\lambda^{(T)}_{N,M} = \sup_{x \in [-1,1]} \sum_{|n| \leq N} \left | W^{(T)}_{N,M}(\mathbf{e}_n)(x) \right |,
}
where $W^{(T)}_{N,M} = R^{(T)}_{N}  - R^{(T)}_{N} \circ L^{(T)}_{N,M} \circ S_M \circ R^{(T)}_{N}$.
}
Much as before, we may now write
\bes{
\lambda^{(T)}_{N,M} = \lim_{K \rightarrow \infty} \lambda^{(T)}_{N,M,K},\qquad \lambda^{(T)}_{N,M,K} = \max_{k=1,\ldots,K_T} \sum_{|n| \leq N} \left | W^{(T)}_{N,M}(\mathbf{e}_n)(x_k) \right |.
}
where the latter can once more be computed efficiently using FFTs.

\begin{figure}
\begin{center}
$\begin{array}{ccc}
\includegraphics[width=5.00cm]{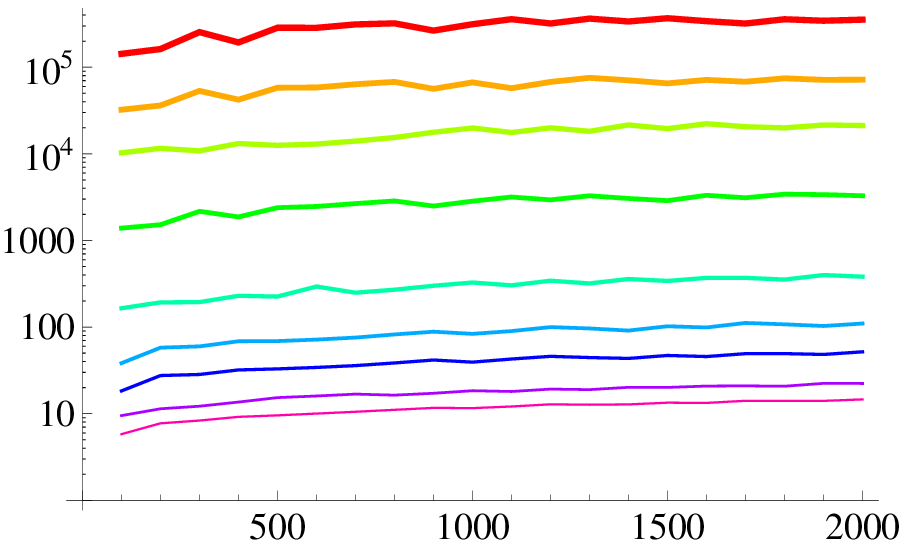}  &  \includegraphics[width=5.00cm]{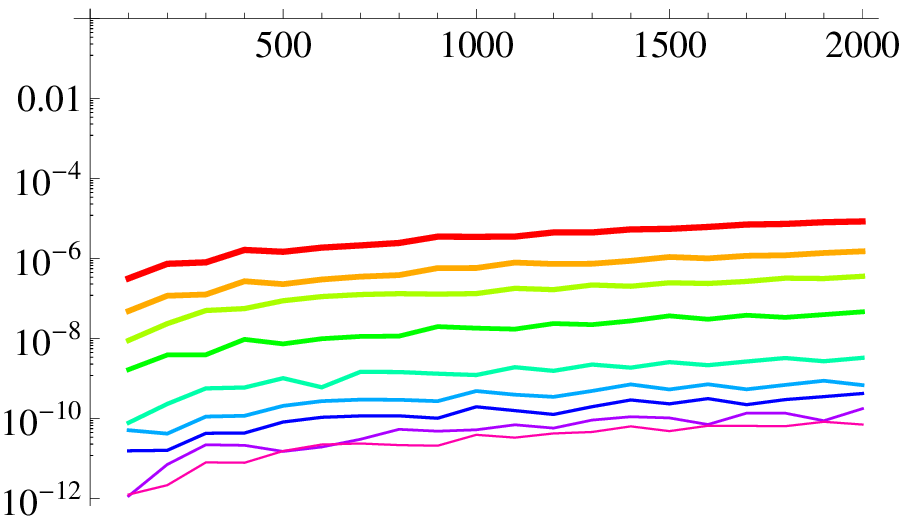} &\includegraphics[width=5.00cm]{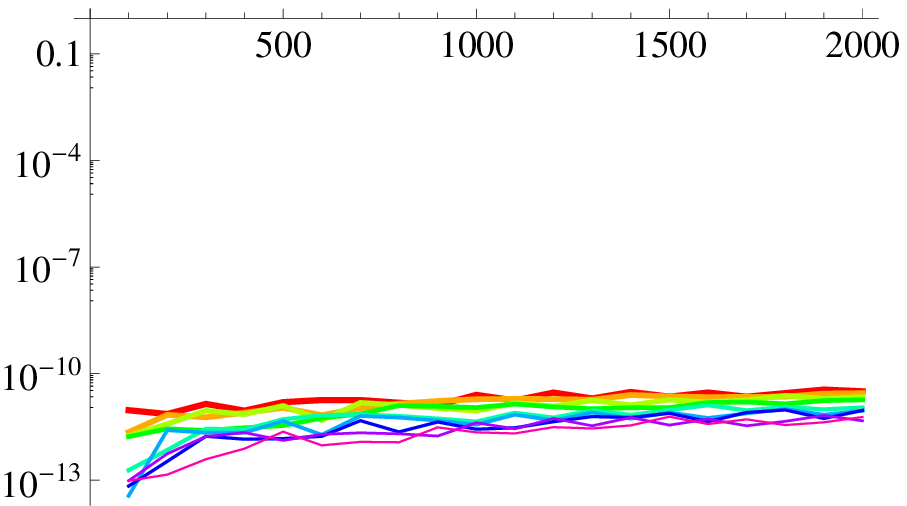}  
\\ 
 \includegraphics[width=5.00cm]{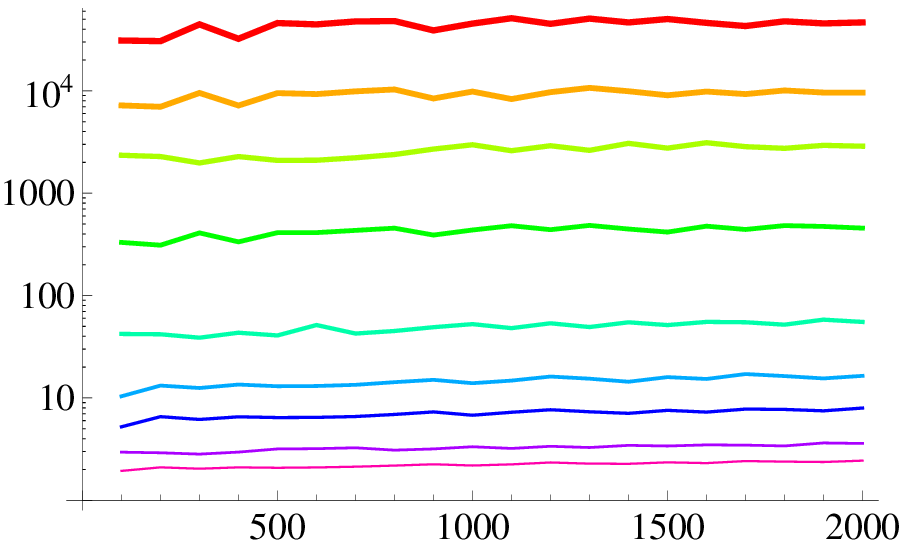}  &  \includegraphics[width=5.00cm]{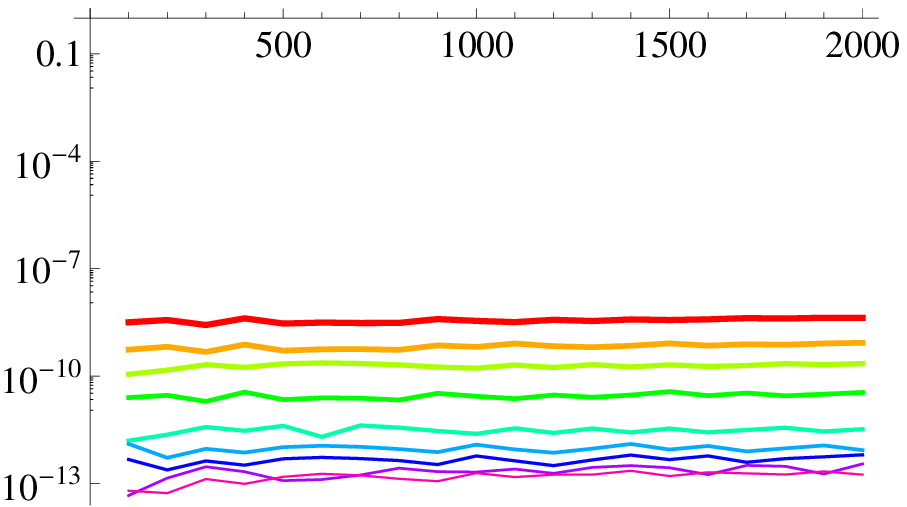} &  \includegraphics[width=5.00cm]{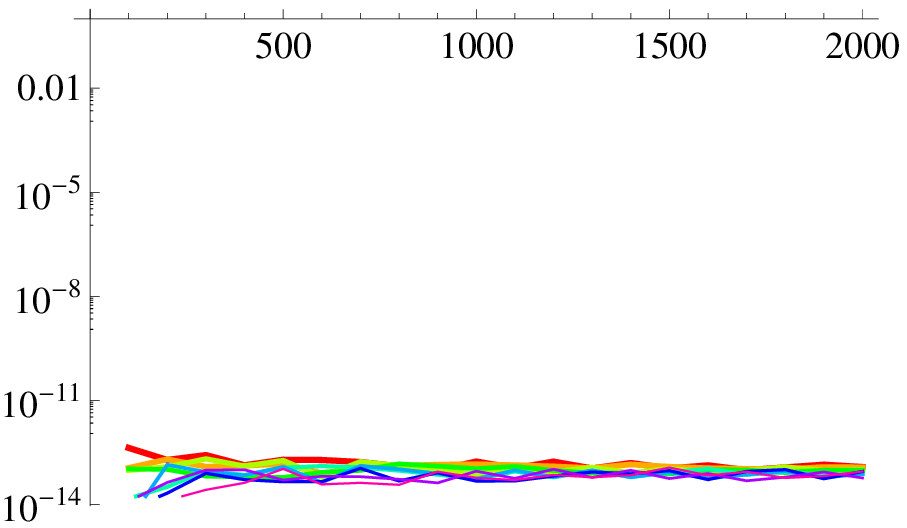}  \\ 
\end{array}$
\caption{Top row: the quantities $\kappa^{(T,\epsilon)}_{M/\eta ,M}$ (left), $\lambda^{(T,\epsilon)}_{M/\eta ,M}$ (middle) and $\mu^{(T,\epsilon)}_{M/\eta , M} : = \lambda^{(T)}_{M/\eta,M} / \kappa^{(T)}_{M/\eta ,M}$ (right) against $M$ for $\eta = 1,1.125,1.25,1.5,2,2.5,3,4,5$ (thickest to thinnest), $\epsilon = 10^{-13}$ and $T=2$.  Bottom row: the same quantities scaled by $\log M$, $M$ and $M / \log M$ respectively.} \label{f:T_2_kappa_lambda}
\end{center}
\end{figure}

\begin{figure}
\begin{center}
$\begin{array}{ccc}
\includegraphics[width=5.00cm]{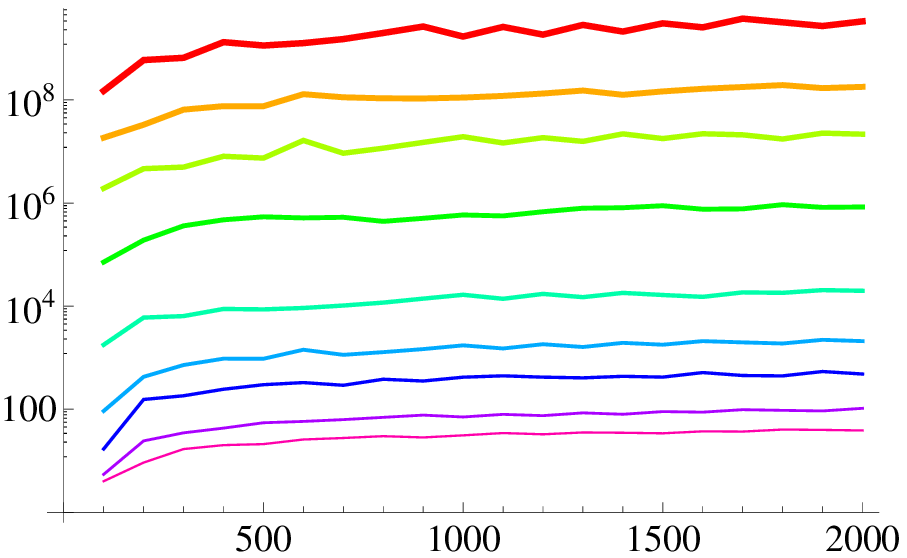}  &  \includegraphics[width=5.00cm]{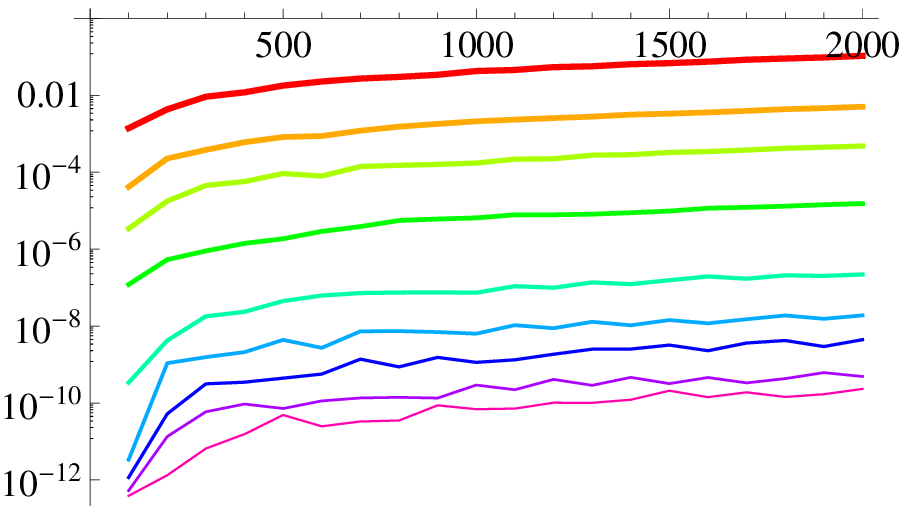} & \includegraphics[width=5.00cm]{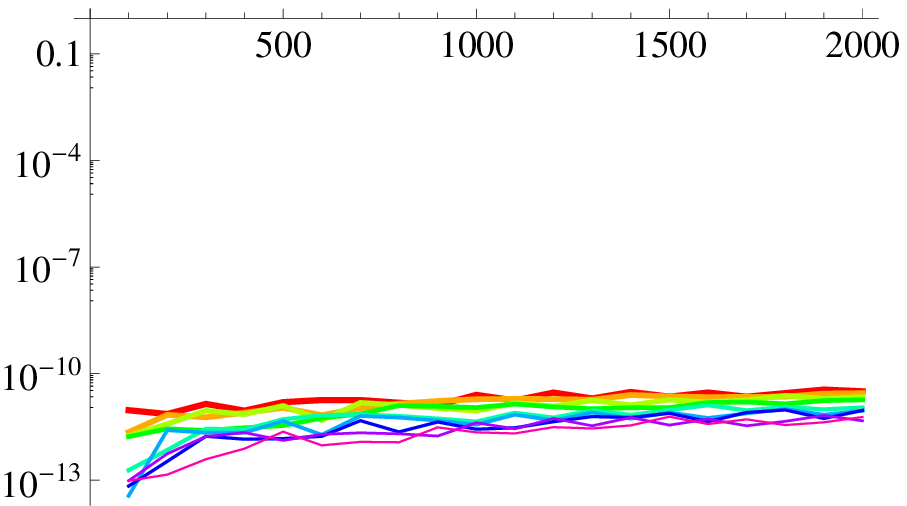}  
\\ 
 \includegraphics[width=5.00cm]{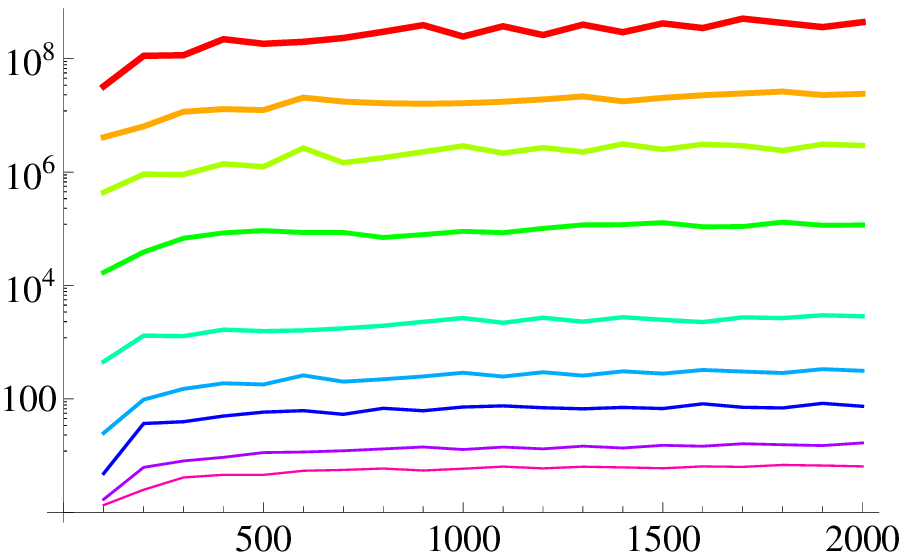}  &  \includegraphics[width=5.00cm]{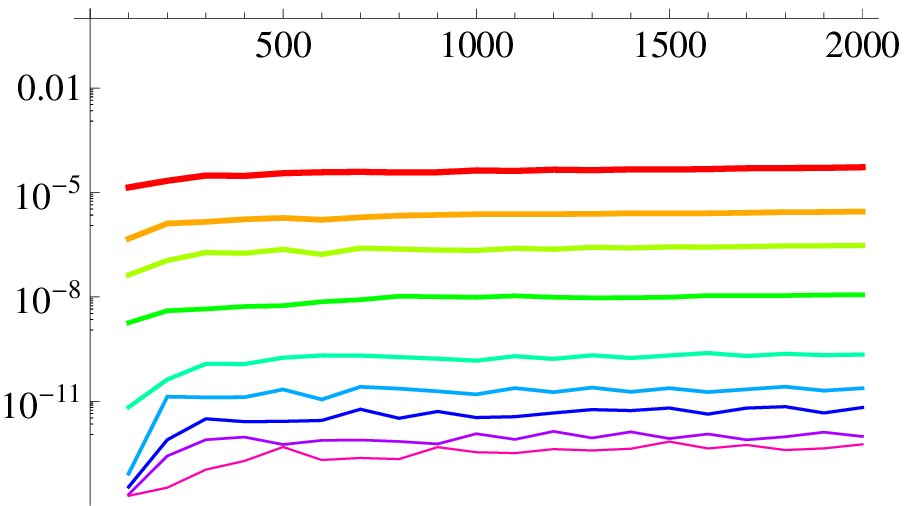} &  \includegraphics[width=5.00cm]{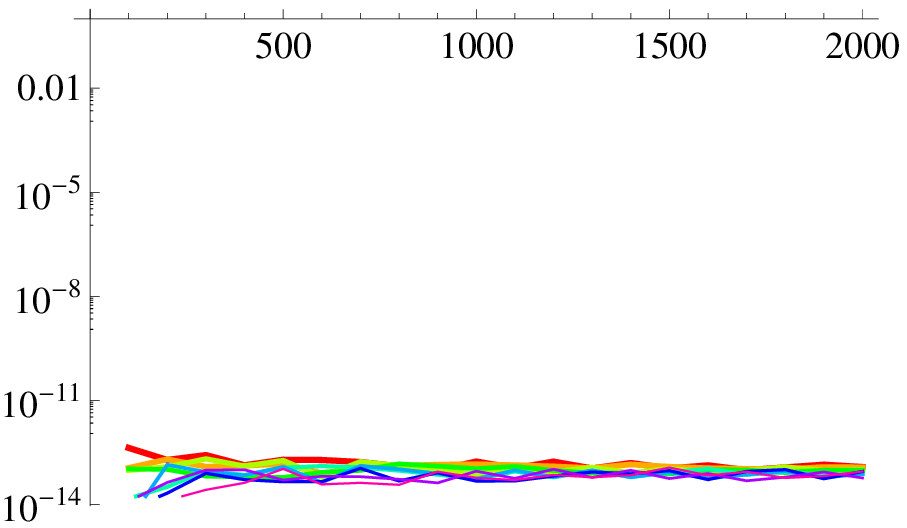}  \\ 
\end{array}$
\caption{Top row: the quantities $\kappa^{(T,\epsilon)}_{M/\eta ,M}$ (left), $\lambda^{(T,\epsilon)}_{M/\eta ,M}$ (middle) and $ \lambda^{(T)}_{M/\eta,M} / \kappa^{(T)}_{M/\eta ,M}$ (right) against $M$ for $\eta = 1,1.125,1.25,1.5,2,2.5,3,4,5$ (thickest to thinnest), $\epsilon = 10^{-13}$ and $T=1.25$.  Bottom row: the same quantities scaled by $\log M$, $M$ and $M / \log M$ respectively.} \label{f:T_125_kappa_lambda}
\end{center}
\end{figure}

Having demonstrated how to compute $\kappa^{(T)}_{N,M}$ and $\lambda^{(T)}_{N,M}$, in Figures \ref{f:T_2_kappa_lambda} and \ref{f:T_125_kappa_lambda} we present the result of such computations for the case where the FE is computed using an SVD with tolerance $\epsilon = 10^{-13}$.  The results confirm the scaling \R{kappa_lambda_growth} for these quantities.  Moreover, these results also show that the quantity
\be{
\label{mu_behaviour}
 \mu^{(T)}_{M/\eta ,M} =  \left (\frac{\lambda^{(T)}_{M/\eta,M}}{\kappa^{(T)}_{M/\eta ,M}} \right ) \left ( \frac{\log M}{M} \right ) \leq \mu,\quad \forall M \gg 1,
}
where $\mu$ is roughly $10^{-13}$ in magnitude, regardless of the choice of $T$ and $\eta$.

Note that in both cases a larger oversampling ratio $\eta$ leads to a smaller condition number $\kappa^{(T,\epsilon)}_{M/\eta,M}$ and numerical defect constant $\lambda^{(T,\epsilon)}_{M/\eta,M}$.  Moreover, the larger value of $T$, in this case, $T=2$, has a smaller condition number for the same oversampling value than the smaller value $T=1.25$.  The purpose of the next section is to investigate the exact nature of these relative scalings.

\rem{
Several previous papers have investigated quantities similar to $\kappa^{(T)}_{N,M}$ and $\lambda^{(T)}_{N,M}$.  In \cite{FEStability} and \cite{LyonFESVD}, quantities based on the exact singular values and vectors of the matrix $A^{(T)}$ were investigated using high-precision numerical computations.  The approach we take above differs from these studies in two aspects.  First, $\kappa^{(T)}_{N,M}$ and $\lambda^{(T)}_{N,M}$ can be formulated for any numerical solver, not just truncated SVDs.  Second, when truncated SVDs are used, they incorporate the numerical errors in the calculation of the singular values and vectors.  Since the matrix $A^{(T)}$ is ill-conditioned, these errors cannot be assumed to be insignificant.
}

\section{Numerical investigation}
We now suppose that the FE is computed using an SVD as described in \S \ref{ss:FE_equi}.  Our aim is to examine the behaviour of $\kappa^{(T,\epsilon)}_{N,M}$, and later $\lambda^{(T,\epsilon)}_{N,M}$, with respect to $T$ and the \textit{oversampling} ratio $\eta = M/N$, and how this affects the accuracy of the corresponding FE approximation.

\begin{figure}
\begin{center}
$\begin{array}{ccc}
\includegraphics[width=8.00cm]{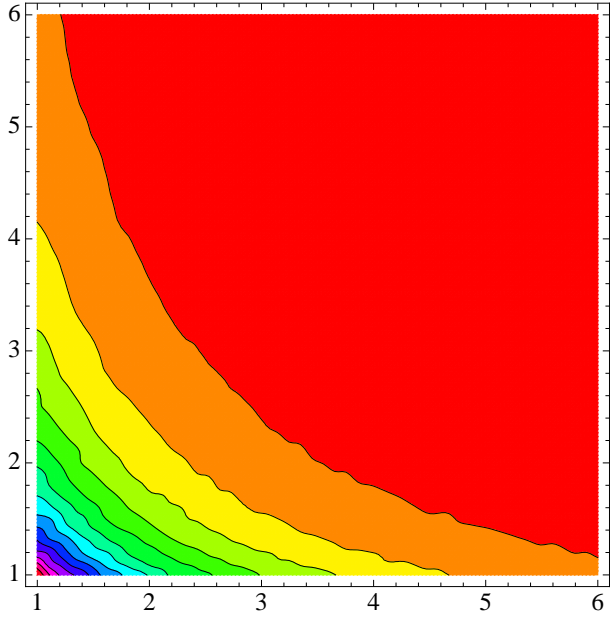}
\end{array}$
\caption{\small Contour plot of $\kappa^{(T,\epsilon)}_{N,\eta N}$ against $1 < T \leq 6$ (horizontal axis) and $1 \leq \eta \leq 6$ (vertical axis) for $N=250$ and $\epsilon = 10^{-13}$.} \label{f:ContourPlot}
\end{center}
\end{figure}

In Figure \ref{f:ContourPlot} we give a contour plot of $\kappa^{(T,\epsilon)}_{N,\eta N}$ as function of $T$ and $\eta$.  As is evident, increasing either $\eta$ or $T$ leads to a smaller condition number.  This suggests that in practice, a balance must be struck between $\eta$ and $T$ so as to get a good condition number whilst retaining a good approximation properties (recall that larger values of $T$ possess worse resolution power \cite{FEStability}, whereas larger values of $\eta$ yield worse approximations for a fixed budget of $M$ equispaced data points).  

Before investigating this interplay further, let us briefly explain why the condition number behaves in this way.  Clearly, increasing $\eta$ results in a smaller value of $N = M/\eta$, and therefore an approximation space $\cG^{(T)}_{N}$ of smaller dimension.  The condition number, the maximum taken over this space, therefore decreases.  Conversely, when $T$ is decreased, this means that the Fourier basis is defined over a smaller domain.  But the Fourier series $\phi = R^{(T)}_{N} \circ L^{(T)}_{N,M}(\mathbf{b})$ must fit the nonperiodic data $\mathbf{b}$ in a least-squares sense on the original domain $[-1,1]$ and must be periodic on the extended domain.  To do this, $\phi$ will be required to take increasingly larger values between data points as $T \rightarrow 1^+$, giving it a bigger uniform norm in comparison to the data norm $| \mathbf{b} |_{\infty}$.

\subsection{Setup}
In order to make a comparison, for each different value of $T$ used we shall choose the ratio $\eta = M/N$ in such a way that the condition number is the same.  Specifically, let $\kappa^* >1$ be fixed.  Then for each $T$ and each $M$ in some specified range, we numerically compute the maximum $N$ such that the condition number is no more than $\kappa^* \log M$.  In other words, we compute the function
\be{
\label{Theta_def}
\Theta^{(T)}(M;\kappa^*) = \max \left \{ N : \kappa^{(T)}_{N,M} \leq \kappa^* \log M \right \},\quad M \in \bbN.
}
Note that we allow $\log M$ factor here since $\kappa^{(T)}_{M/\eta,M}$ grows like $\log M$ as $M \rightarrow \infty$.  With this scaling, $\Theta^{(T)}(M;\kappa^*)$ will be linear in $M$ (see later).

The function \R{Theta_def} can be computed numerically for each $M$.  This follows from the fact that $\kappa^{(T)}_{N,M}$ can be computed (see \S \ref{ss:computing}).  Note also that $\Theta^{(T)}(M;\kappa^*)$ can be computed for any particular numerical solver used to solve the least squares \R{LNM_def}, and thus allows a comparison between different methods.  We consider this further in \S \ref{ss:other_solvers}.

Having computed \R{Theta_def} for each value of $T = T_1,\ldots,T_r$ and some range of $M$, we next use these values to compare approximation properties of the corresponding equispaced FEs
\bes{
F^{(T_j)}_{N,M}(\cdot),\quad \mbox{where $N = \Theta^{(T_j)}(M;\kappa^*)$},\qquad j=1,\ldots,r.
}
When doing this, we shall consider a suite of different test functions, described further below.  Observe that \R{Theta_def} determines the largest value of $N$ for which the condition number is at most $\kappa^* \log M$.  Thus, setting $N = \Theta^{(T)}(M;\kappa^*)$ when computing the FE $F^{(T)}_{N,M}(f)$ ensures the best approximation properties for each value of $T$ (since $N$ is maximal) whilst retaining the same condition number for the different choices $T=T_1,\ldots,T_r$.  Thus a comparison between these different values of $T$ can be made, using the condition number as the common fixed point.

\begin{figure}
\begin{center}
$\begin{array}{ccc}
\includegraphics[width=5.00cm]{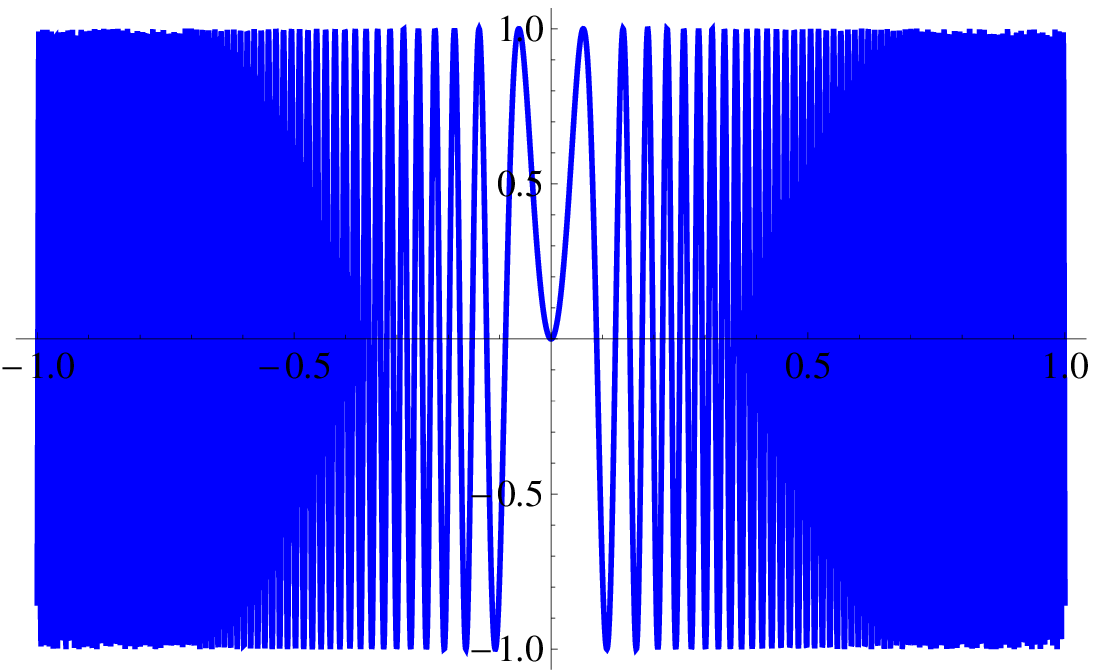} & \includegraphics[width=5.20cm]{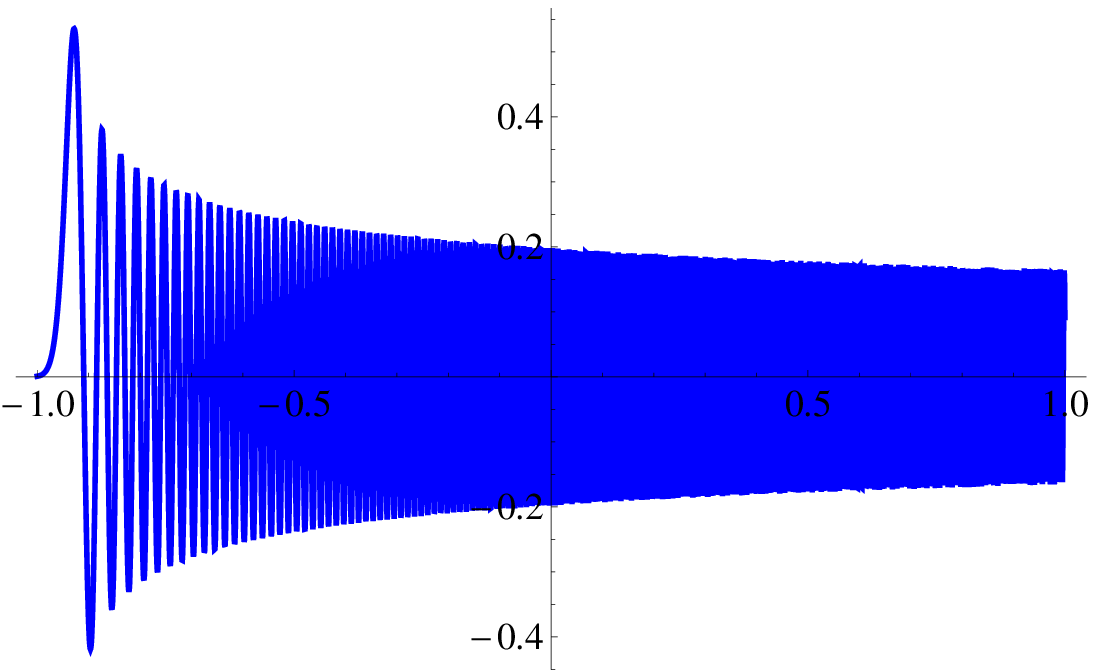} & \includegraphics[width=5.20cm]{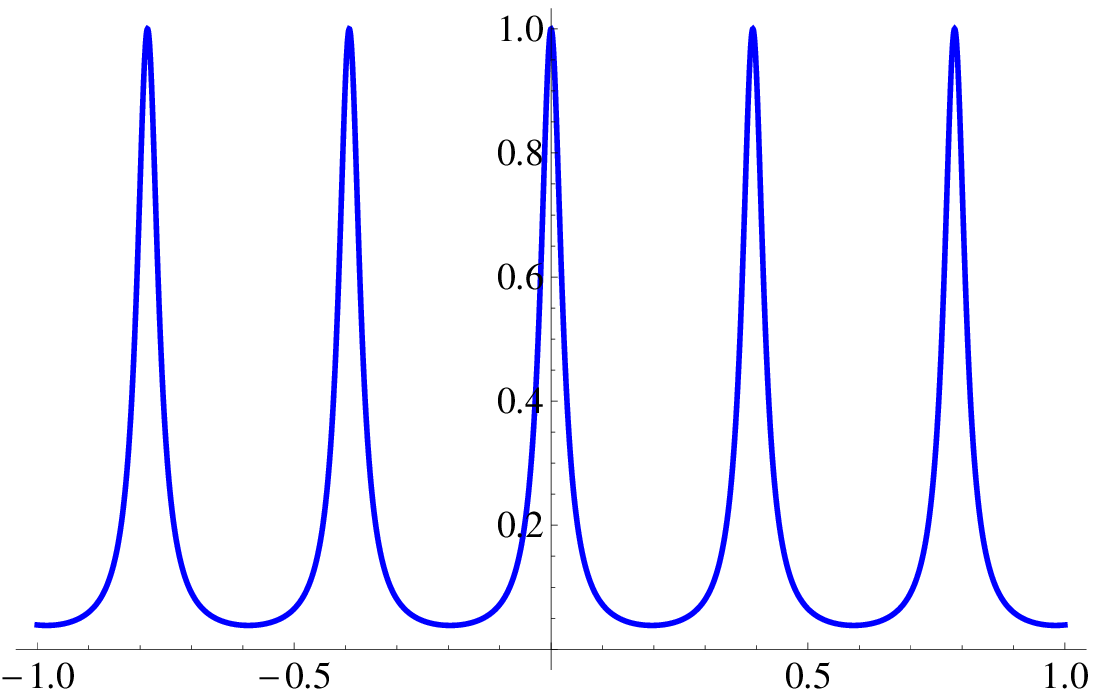} \\
\includegraphics[width=5.00cm]{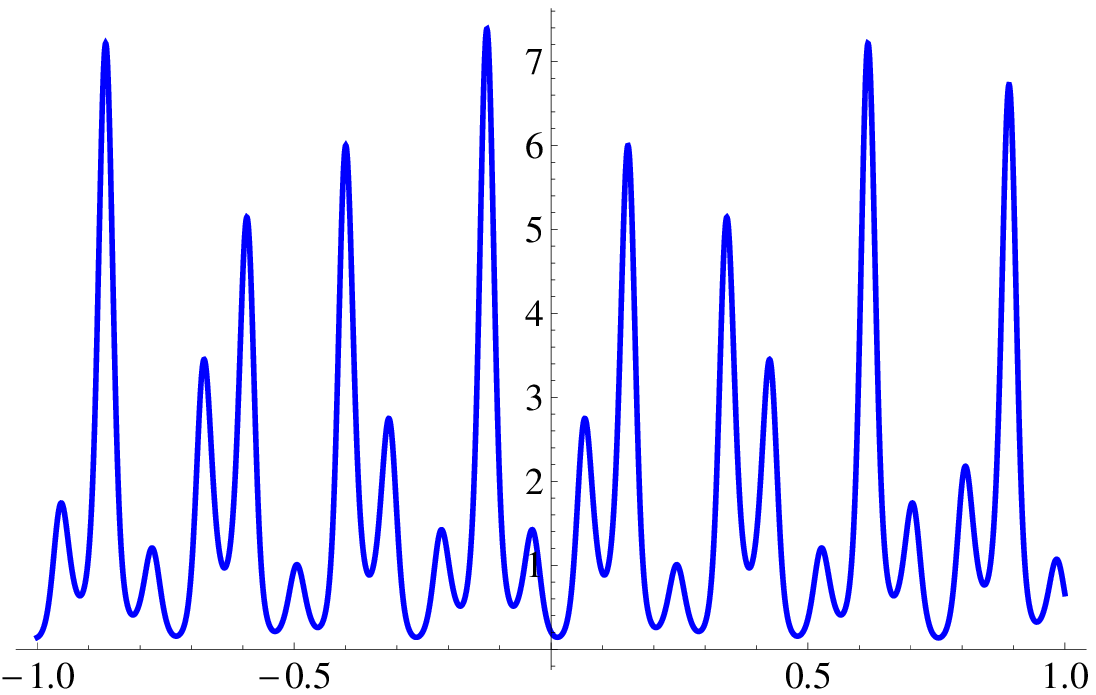} & \includegraphics[width=5.20cm]{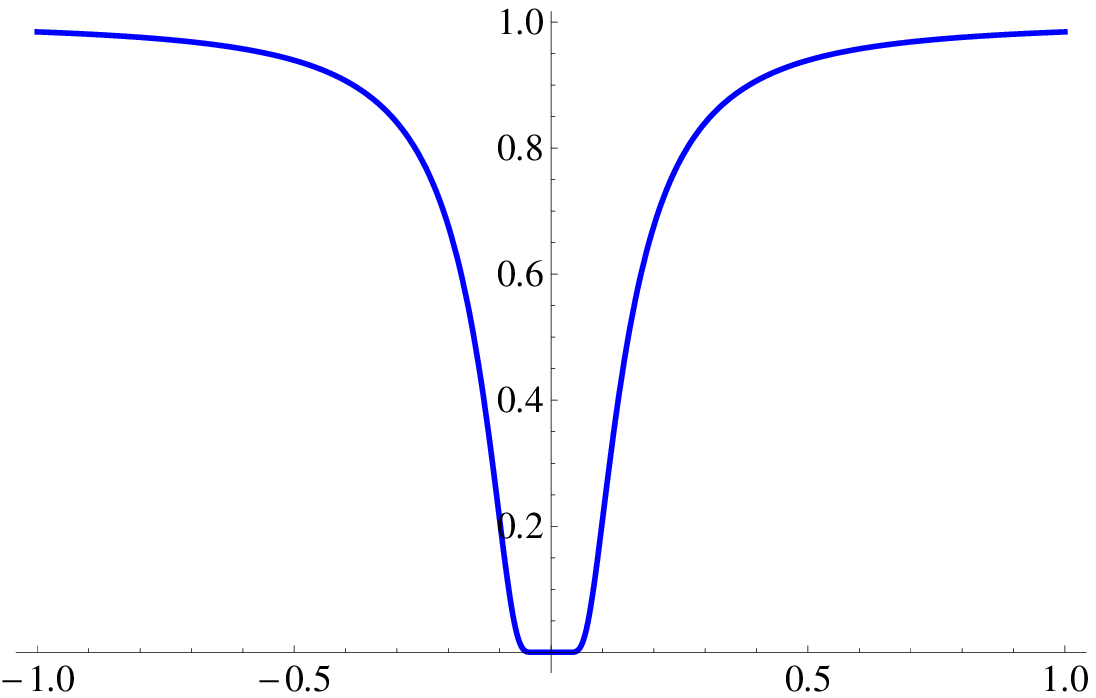} & \includegraphics[width=5.20cm]{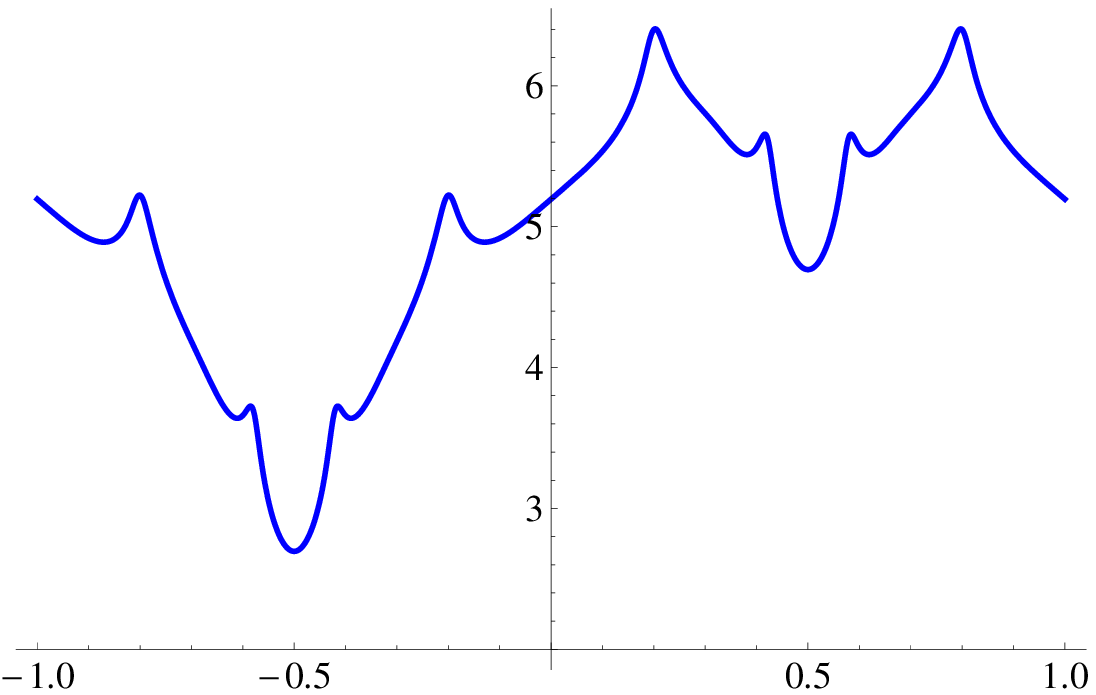}
\end{array}$
\caption{Top row: the functions $f_2$, $f_3$ and $f_6$.  Bottom row: the functions $f_7$, $f_8$ and $f_9$.} \label{f:functions}
\end{center}
\end{figure}

We now require appropriate test functions.  Our first three functions are as follows: 
\eas{
f_1(x) = \E^{230 \sqrt{2} \I \pi x}\qquad f_2(x) = \sin(400 x^2)\qquad f_3(x) = \mathrm{Ai}(-66-70 x)
}
These functions all exhibit oscillations, which make them challenging to approximate from equispaced data.  Plots of $f_2$ and $f_3$ are given in Figure \ref{f:functions}.  Our next collection of test functions feature singularities in the complex plane near $[-1,1]$, again making them difficult to approximate:
\eas{
f_4(x) = \frac{1}{1+1500 x^2}\qquad f_5(x) = \frac{1}{60-59 x}\qquad f_6(x) = \frac{1}{1+25 \sin^2 8 x}.
}
A plot of $f_6$ is shown in Figure \ref{f:functions}.  Our final collection of functions, also displayed in Figure \ref{f:functions}, is 
\bes{
\quad f_7(x) =  \E^{\sin (21.6 \pi x - 10.8 \pi )-\cos 8 \pi x}\qquad f_8(x) = \E^{-1/(8 x)^2}\qquad f_9(x) = s(x).
}
Note that $f_7$ is often used in testing algorithms for recovering functions to high accuracy from equispaced data, and $f_8$ is made challenging by its lack of analyticity and the flat region near $x=0$.  The function $f_9$ is similar to that introduced in \cite{TrefethenFunctions}.  It is obtained by the following iteration:
\bull{
\item $s(x) = \sin \pi x$, $f(x) = s(x)$
\item For $j=1,2,\ldots,10$, $s(x) = 3/4 (1-2 s(x)^4)$, $f(x) = f(x) + s(x)$. 
}

\subsection{Numerical results}

\begin{figure}
\begin{center}
$\begin{array}{ccc}
\includegraphics[width=5.00cm]{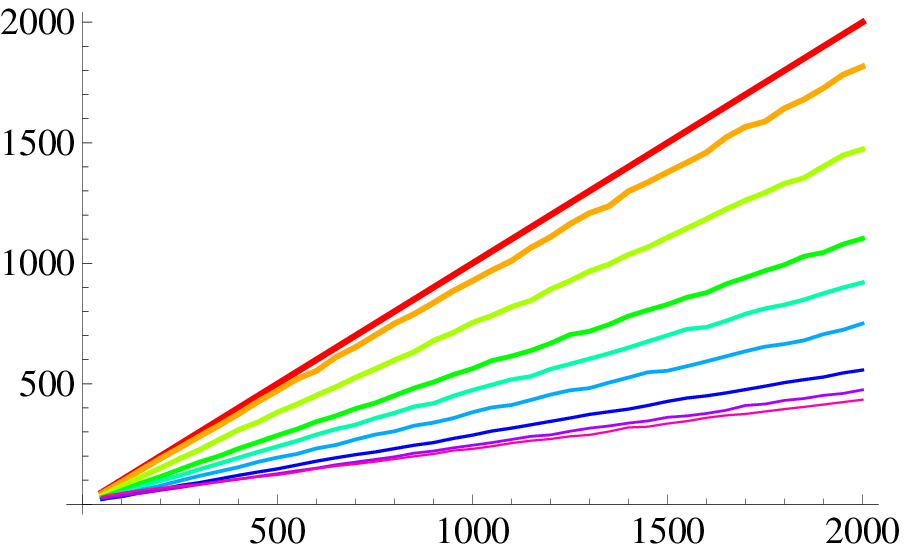} &
\includegraphics[width=5.00cm]{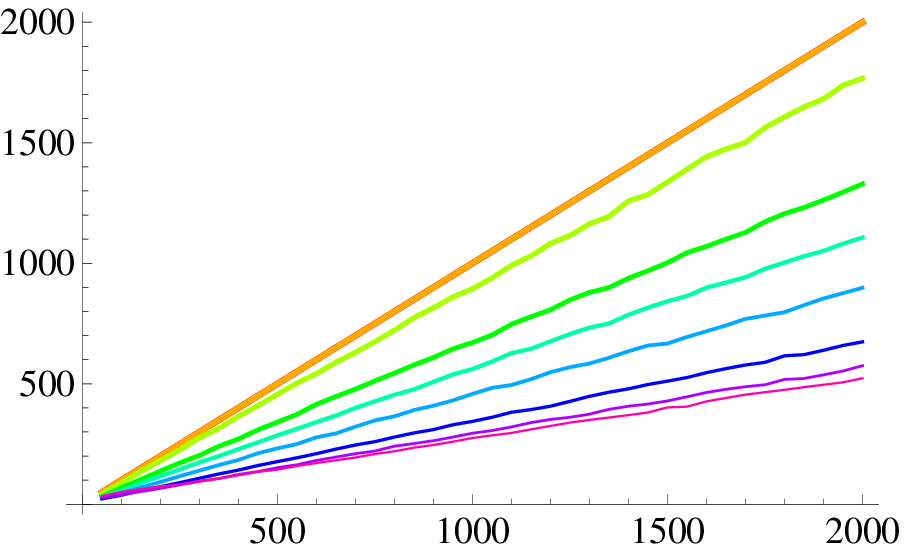} &
\includegraphics[width=5.00cm]{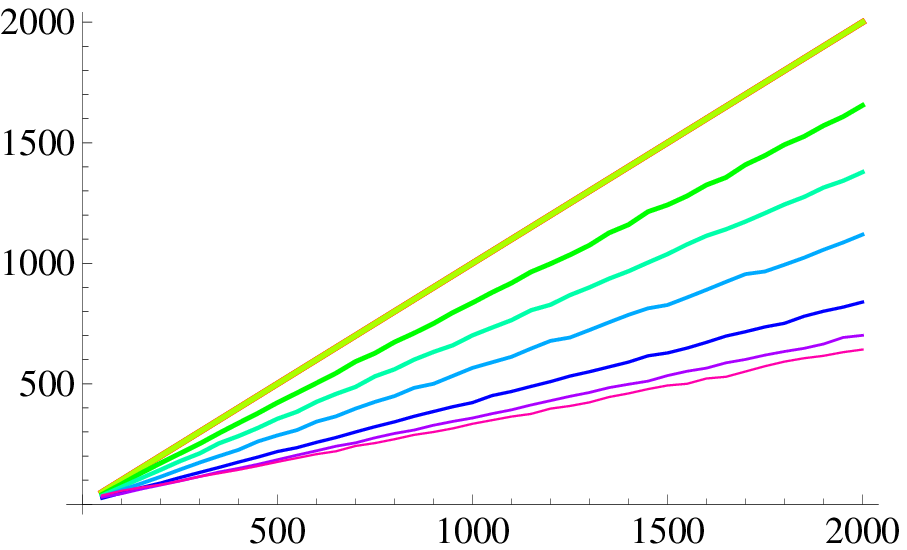} 
\\ 
\includegraphics[width=5.00cm]{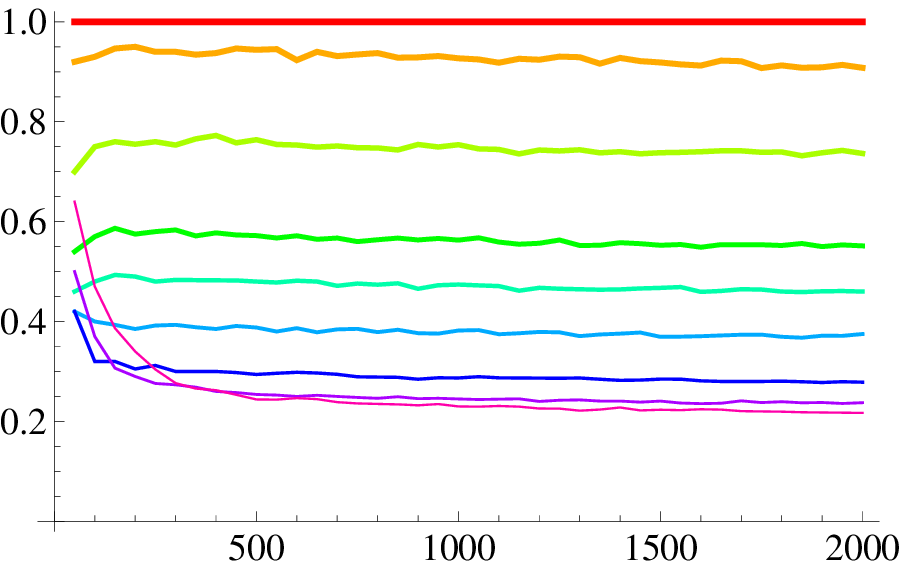} &
\includegraphics[width=5.00cm]{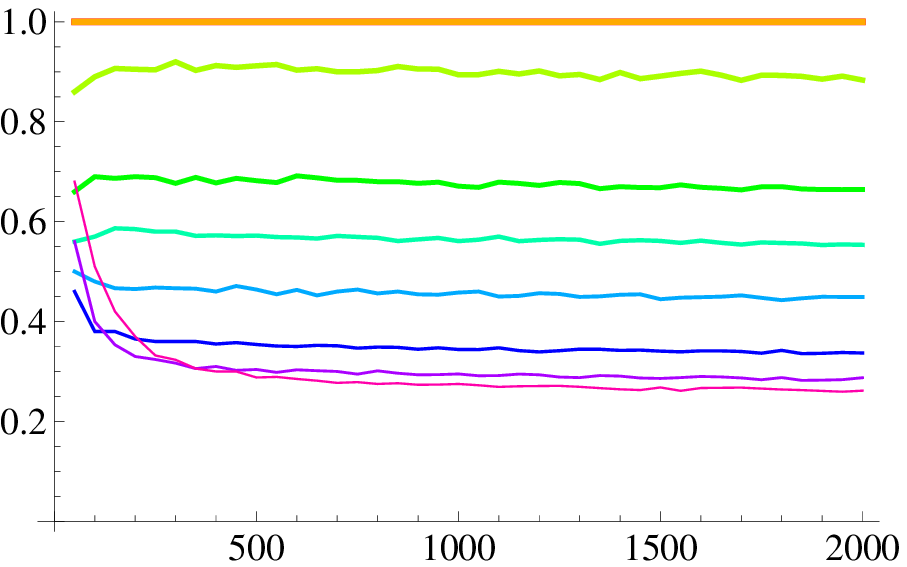} &
\includegraphics[width=5.00cm]{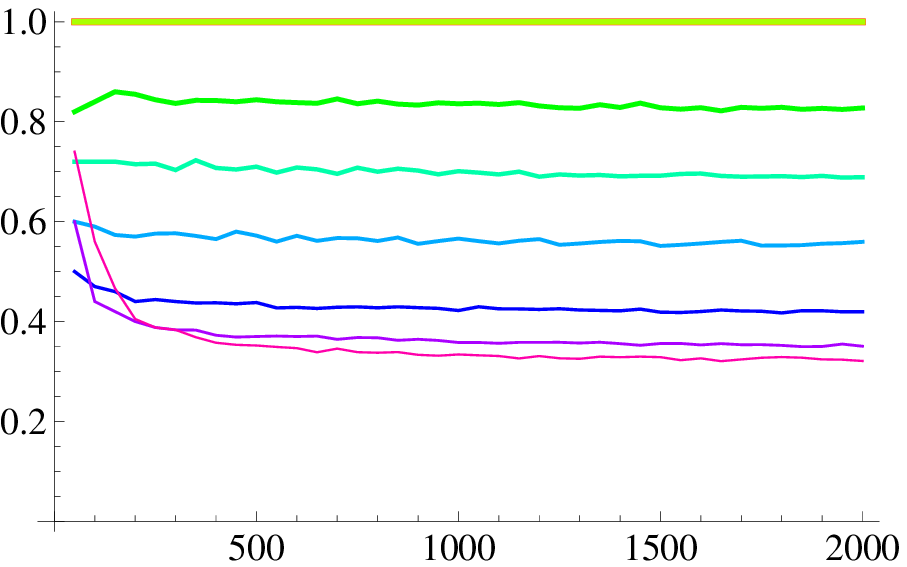} 
\end{array}$
\caption{\small Plots of $\Theta^{(T,\epsilon)}(M;\kappa^*)$ (top) and $\Theta^{(T,\epsilon)}(M;\kappa^*)/M$ (bottom) against $M$ for $\kappa^* = 10$ (left), $\kappa^* =25$ (middle) and $\kappa^* = 100$ (right) using $\epsilon = 10^{-13}$.  The values of $T$ used (in order of increasing thickness) were $T=1.125,1.25,1.5,2,2.5,3.0,4.0,5.0,6.0$.  Note that in the middle plots the $T=5.0$ line is identical to the $T=6.0$ line, and for the right plots the $T=4.0$, $T=5.0$ and $T=6.0$ lines are identical.} \label{f:ThetaPlot}
\end{center}
\end{figure}

In Figure \ref{f:ThetaPlot} we plot the function $\Theta^{(T)}(M;\kappa^*)$ against $M$ for various values of $T$ and $\kappa^*$.  Note that this function is approximately linear in $M$.  Moreover, its gradient is larger for bigger values of $T$, as expected from the results given in Figure \ref{f:ContourPlot}.  The approximate linear rate of growth of $\Theta^{(T)}(M;\kappa^*)$ is shown in Table \ref{tab:ThetaGrowth}.

\begin{table}
\begin{center}
\begin{tabular}{|c|c|c|c|c|c|c|c|c|c|}
\hline
$T$ & $1.125$ & $1.25$ & $1.5$ & $2.0$ & $2.5$ & $3.0$ & $4.0$ & $5.0$ & $6.0$ \\
\hline
$\kappa^* = 10$ & 0.21 & 0.23 & 0.28 & 0.37 & 0.46 & 0.55 & 0.73 & 0.91 & 1.00  \\
\hline
  $\kappa^* = 25$ & 0.25 & 0.28 & 0.33 & 0.45 & 0.55 & 0.66 & 0.89 & 1.00 & 1.00 \\
 \hline
$\kappa^* = 100$ &  0.31 & 0.35 & 0.42 & 0.55 & 0.69 & 0.82 & 1.00 & 1.00 & 1.00 \\
\hline
\end{tabular}
\caption{The approximate linear scaling of $\Theta^{(T)}(M;\kappa^*)$ with $M$.  Values were computed using linear regression on the data obtained in Figure \ref{f:ThetaPlot}.}\label{tab:ThetaGrowth}
\end{center}
\end{table}

\begin{figure}
\begin{center}
$\begin{array}{ccc}
\includegraphics[width=5.00cm]{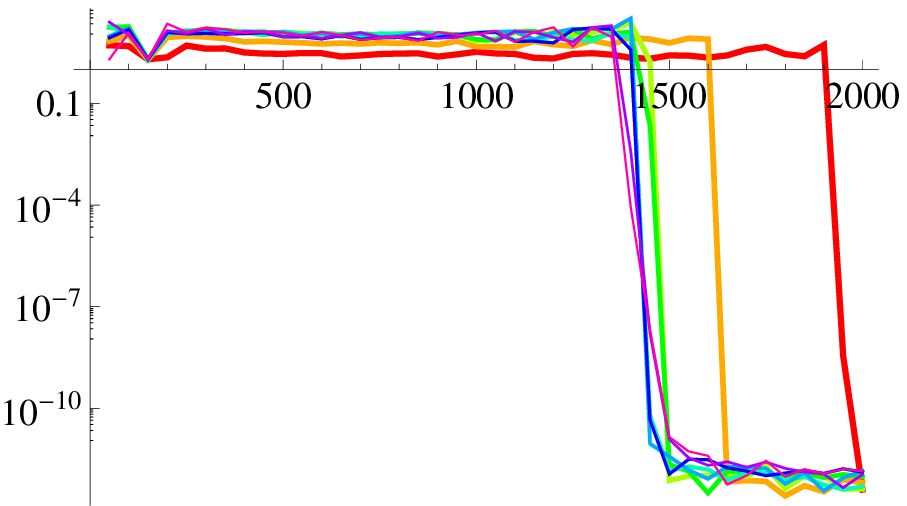}  &
\includegraphics[width=5.00cm]{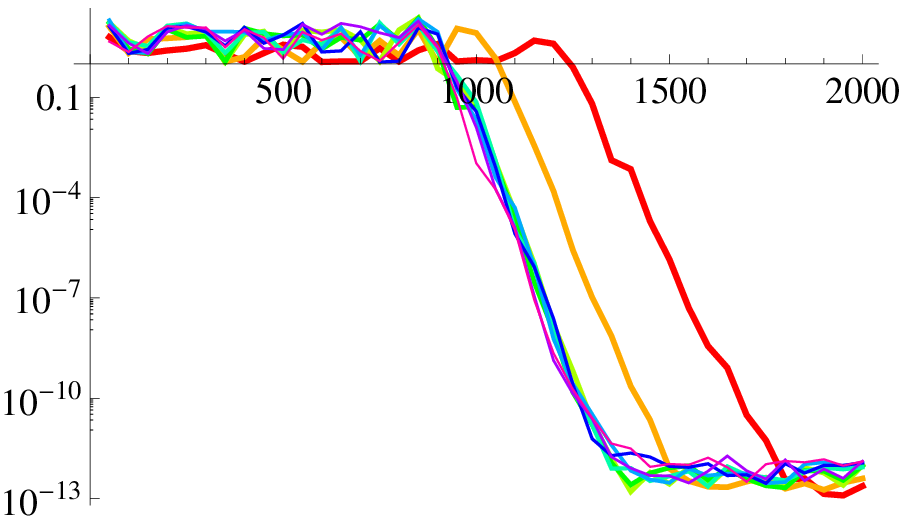}  &
\includegraphics[width=5.00cm]{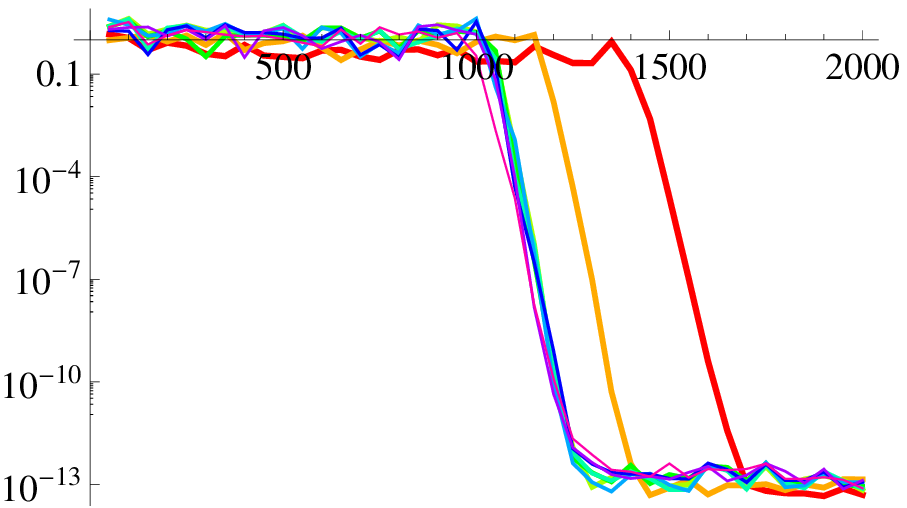}  \\
\includegraphics[width=5.00cm]{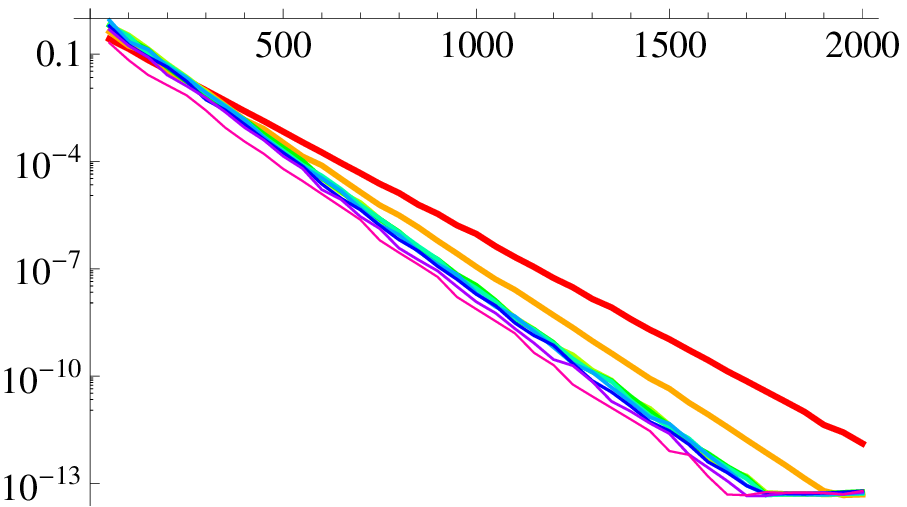}  &
\includegraphics[width=5.00cm]{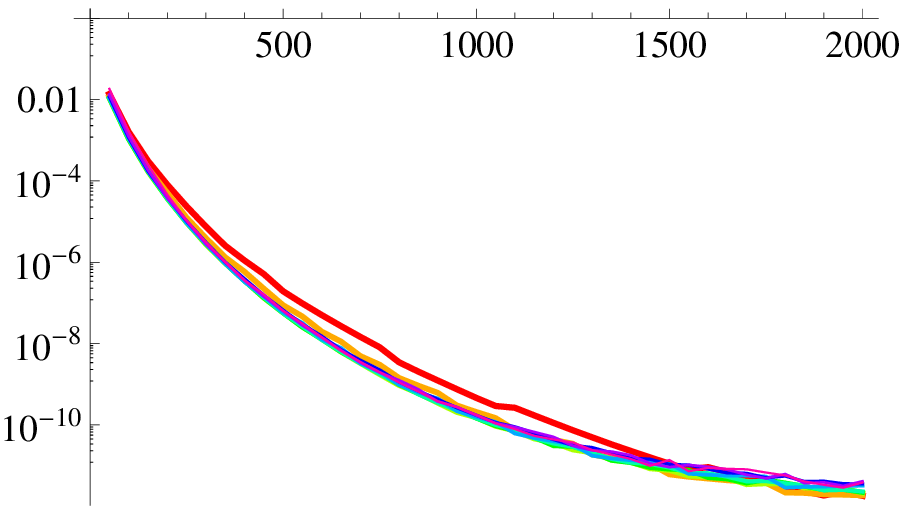}  &
\includegraphics[width=5.00cm]{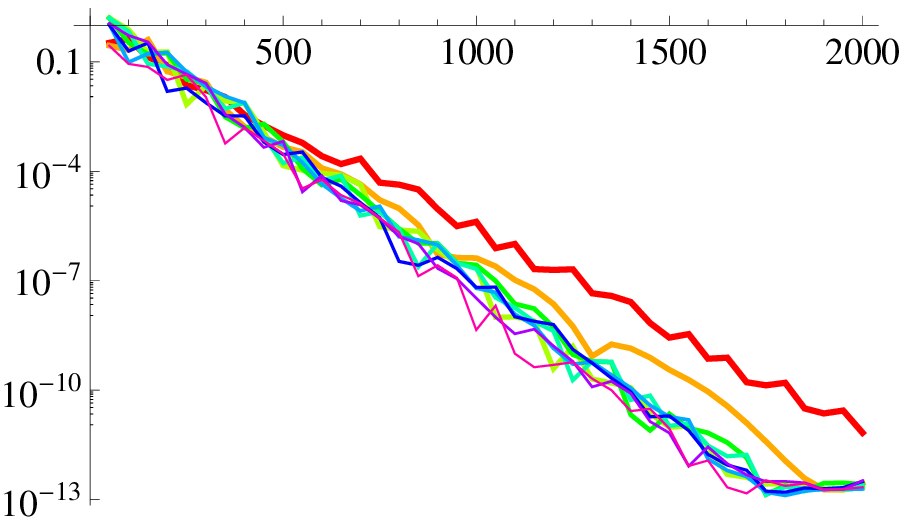}  \\
\includegraphics[width=5.00cm]{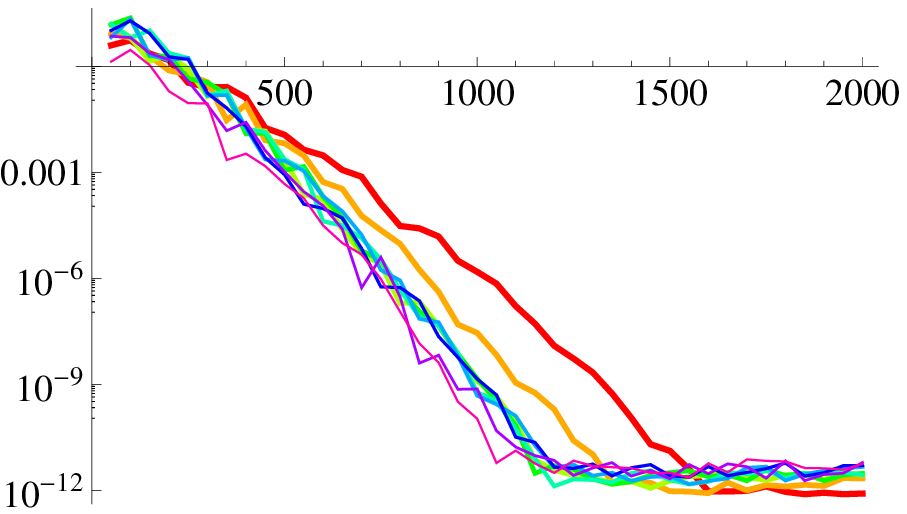}  &
\includegraphics[width=5.00cm]{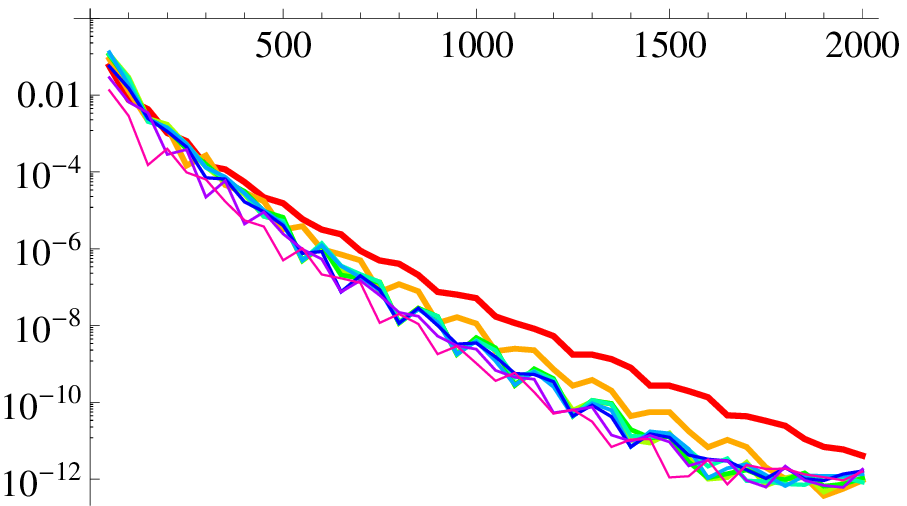}  &
\includegraphics[width=5.00cm]{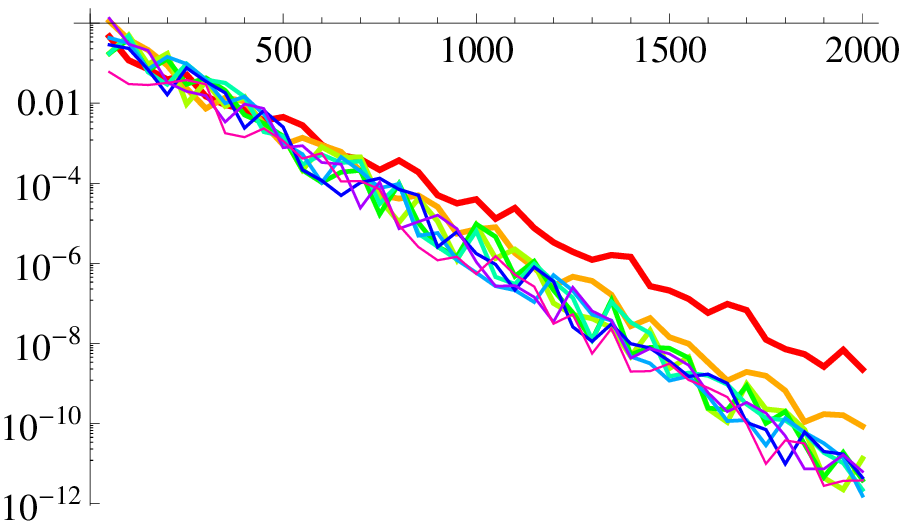}  \\
\end{array}$
\caption{\small Approximation errors for the functions $f_i$, $i=1,\ldots,9$ using the values $\Theta^{(T,\epsilon)}(M;\kappa^*)$ computed in Figure \ref{f:ThetaPlot} for $\kappa^* = 25$ and $\epsilon = 10^{-13}$.} \label{f:FnApp}
\end{center}
\end{figure}

Next, in Figure \ref{f:FnApp} we consider the approximation of the functions $f_1,\ldots,f_9$ using the derived values for $\Theta^{(T)}(M;\kappa^*)$.   These results point towards a surprising phenomenon.  Besides the choices $T=5.0$ and $T=6.0$, all values of $T$ used lead to near-identical approximation errors, regardless of the function considered.  Thus, seemingly the value of $T$, unless taken to be either $5.0$ or $6.0$ in this case, makes little or no difference to the FE approximation. 

\begin{figure}
\begin{center}
$\begin{array}{ccc}
\includegraphics[width=5.00cm]{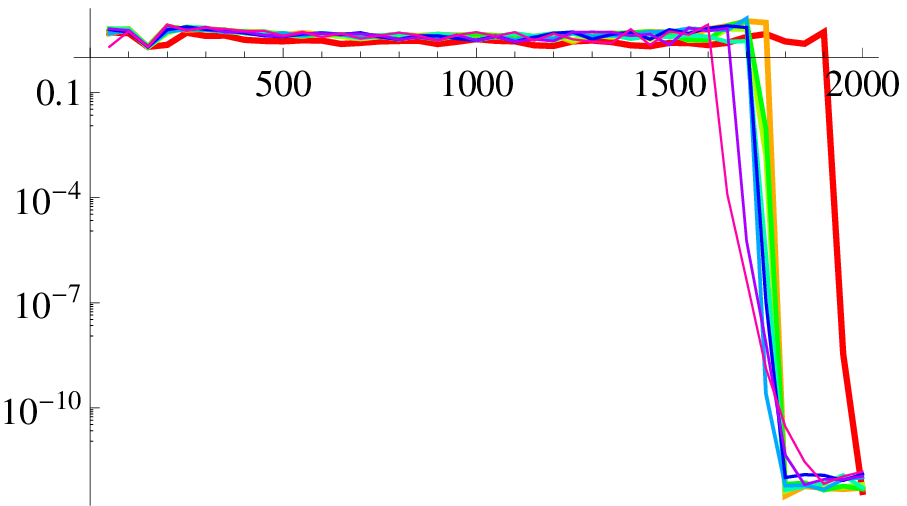}
& 
\includegraphics[width=5.00cm]{f1app_SVD_1e13_kappa25_logscaled} & 
\includegraphics[width=5.20cm]{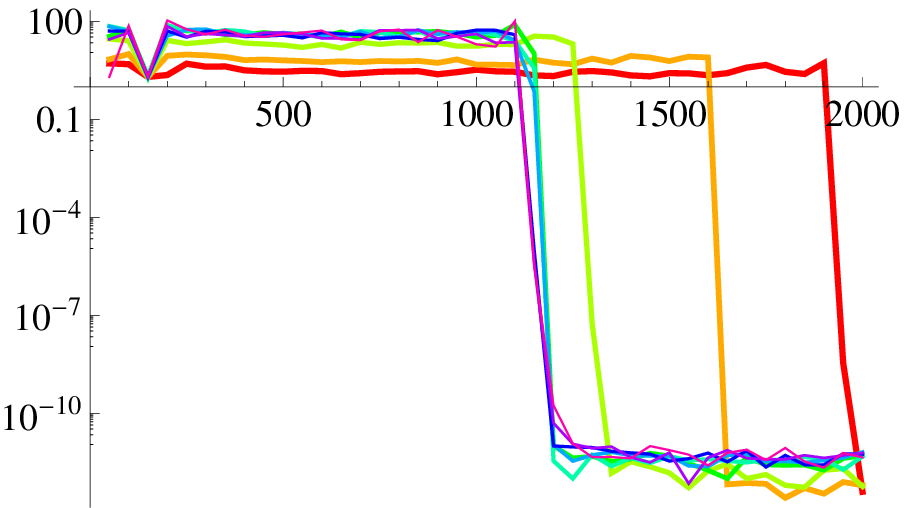}
\end{array}$
\caption{\small Approximation errors for the function $f_1$ using the values $\Theta^{(T)}(M;\kappa^*)$ computed in Figure \ref{f:ThetaPlot} for $\kappa^* = 10$ (left), $\kappa^* = 25$ (middle) and $\kappa^* = 100$ (right).} \label{f:FnApp2}
\end{center}
\end{figure}

Let us explain why the two values $T=5.0$ and $T=6.0$ lead to worse approximations.  This is due to a phenomenon we refer to as \textit{saturation}:

\defn{
For a given $T >1$ and $\kappa^* > 1$, saturation occurs if 
\bes{
\limsup_{M \rightarrow \infty} \kappa^{(T)}_{M,M} / \log M < \kappa^* .
}
}
Saturation means that the corresponding FE is too stable to take advantage of the allowed condition number $\kappa^*$.  In other words, the maximal value of $N$ permitted is limited by the fact that $N \leq M$ in the FE approximation, and not by the condition number constraint $\kappa^{(T)}_{N,M} \leq \kappa^*$.  When saturation occurs, the resulting FE approximation performs worse  in terms of approximation than that corresponding to a value of $T$ for which saturation does not occur.  Figure \ref{f:ThetaPlot} illustrates that the two values of $T$ which give worse approximations in Figure \ref{f:FnApp}, i.e.\ $T=5.0$ and $T=6.0$, do indeed saturate.  This is further demonstrated in Figure \ref{f:FnApp2}, where we consider the approximation of the function $f_1$ using different values of $\kappa^*$.  Figure \ref{f:ThetaPlot} shows that for $\kappa^* =10$ only $T=6.0$ saturates, whereas $T=5.0$ and $T=6.0$ both saturate for $\kappa^* = 25$, and for $\kappa^* = 100$ the values $T=4.0$, $T=5.0$ and $T=6.0$ all saturate.  In Figure \ref{f:FnApp2} the FE approximations with these values of $T$ perform worse than the FEs corresponding to values of $T$ which do not saturate.  Note that when $\kappa^* = 10$ the effect of the saturation for $T=6.0$ has less impact, since $\kappa^{(6,\epsilon)}_{M,M} /\log M \approx 5$ is reasonably close to $\kappa^*$ in this case.  Hence saturation occurs, but to a lesser extent.  Similarly, since $\kappa^{(6)}_{M,M} < \kappa^{(5)}_{M,M} < \kappa^{(4)}_{M,M}$ the effect of the saturation on the FE approximation when $\kappa^* = 100$ is less for $T=4.0$ than it is for $T=5.0$ and $T=6.0$.

With this in hand, we are now able to state the main empirical conclusion of this section: \textit{unless saturation occurs, the choice of $T$ makes little difference to the FE approximation}.  In particular, one may use the value $T=2$ provided it does not saturate for given value of $\kappa^*$.  We note in passing that this phenomenon is of course asymptotic in $M$, and relies on the fact that the functions under consideration require large numbers of equispaced points to be approximated to any accuracy.  For smooth function lacking unpleasant features such as close singularities or oscillations, i.e.\ functions that can be resolved with small $M$, there will be slight discrepancies for different values of $T$.

\begin{figure}
\begin{center}
$\begin{array}{ccc}
\includegraphics[width=5.00cm]{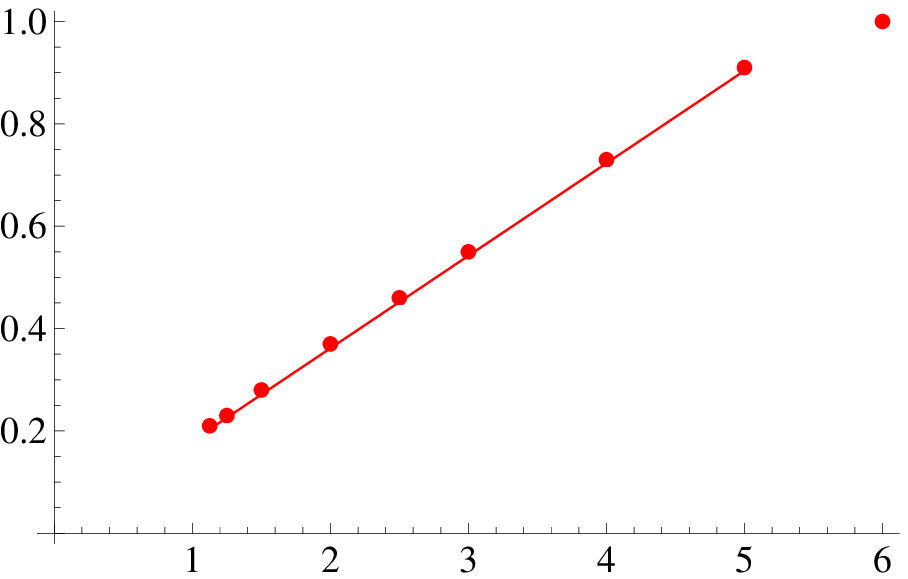}  & 
\includegraphics[width=5.00cm]{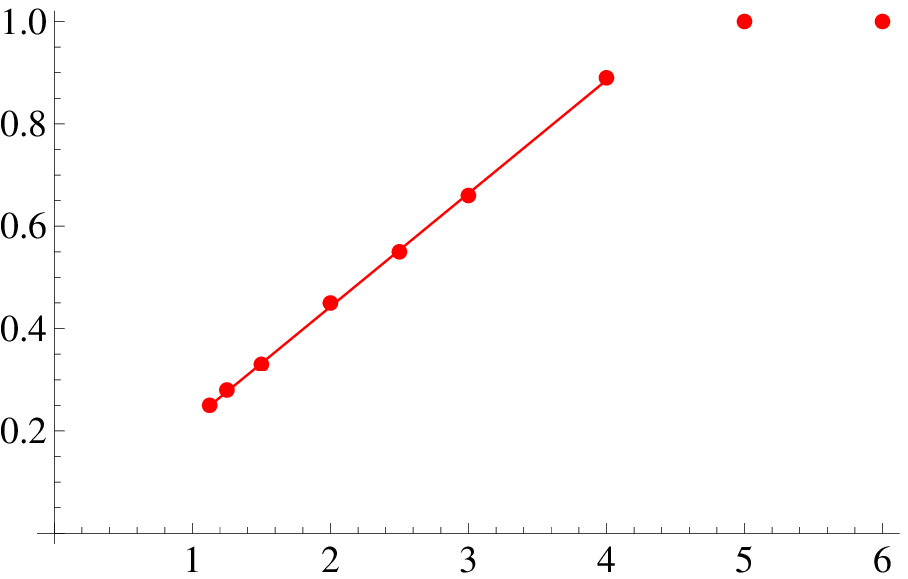}&
\includegraphics[width=5.00cm]{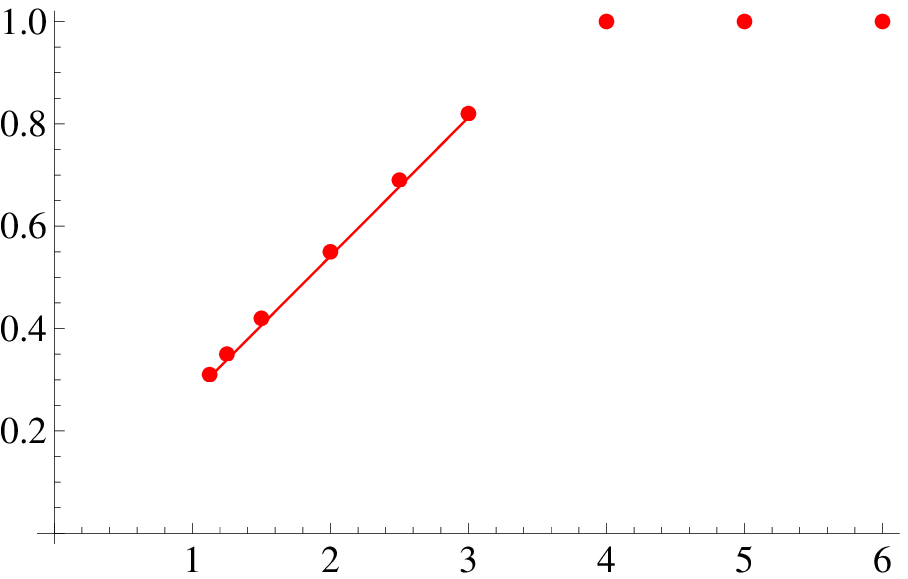} \\
\kappa^* = 10,\ \tau(\kappa^*) \approx 0.18 & \kappa^* = 25,\ \tau(\kappa^*) \approx 0.22 & \kappa^* = 100,\ \tau(\kappa^*) \approx 0.27
\end{array}$
\caption{\small The approximate values $\nu^{(T)}(\kappa^*)$ from the Table \ref{tab:ThetaGrowth} plotted against $T$.  The solid curve is the line with slope $\tau(\kappa^*)$, where $\tau(\kappa^*)$ is computed using linear regression on the data $\{ T , \nu^{(T)}(\kappa^*) \}$ corresponding to those nonsaturating values of $T$.  } \label{f:eta_nu}
\end{center}
\end{figure}

Let us give some explanation for why this conclusion should hold.  Figure \ref{f:ThetaPlot} shows that the function $\Theta^{(T)}(M;\kappa^*)$ is approximately linear in $M$, i.e. $\Theta^{(T)}(M;\kappa^*) \approx \nu^{(T)}(\kappa^*) M$ for large $M$ for some constant $\nu^{(T)}(\kappa^*) >0$.  Explicit values of the quantity $\nu^{(T)}(\kappa^*)$ are given in Table \ref{tab:ThetaGrowth}, and these are plotted against $T$ in Figure \ref{f:eta_nu}.  As is evident from this figure, for each each fixed $\kappa^*$, the quantity $\nu^{(T)}(\kappa^*)$ is approximately linear in $T$ up to the point at which saturation occurs.  Write $\nu^{(T)}(\kappa^*) \approx \tau(\kappa^*) T$ for some $\tau(\kappa^*) > 0$ and now consider the FE approximation $F^{(T)}_{N,M}$ where $N = \tau(\kappa^*) T M$.  By \R{error_numerical_2} and \R{EN_def}, we see that the error is determined up to a mildly growing factor in $M$ by the quantity $E_{N}(f) = E_{\tau(\kappa^*) T M}(f)$.  However, the error bound \R{alg_conv} depends on the ratio $N/T$, which in this case is equal to $\tau(\kappa^*) M$ regardless of the choice of $T$.  Hence the bound \R{alg_conv} is independent of $T$ in the nonsaturating case for this choice of $N$, and therefore we expect the same $T$-independence of the FE approximation whenever $N$ is taken according to $\Theta^{(T)}(M;\kappa^*)$.

\begin{figure}\centering
\includegraphics[width=7.50cm]{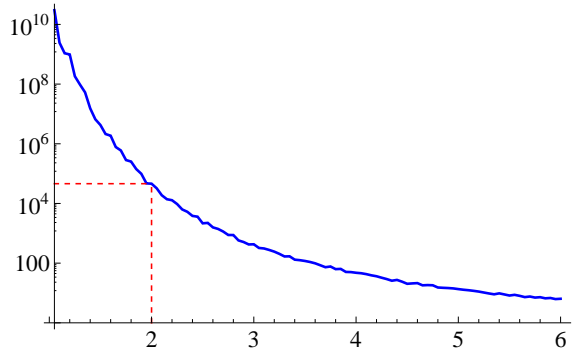}
  \caption{The quantity $\kappa^{(T,\epsilon)}_{M,M} / \log M$ against $1 < T \leq 6$ for $M=500$ and $\epsilon = 10^{-13}$.} \label{f:Saturation_Compute}
\end{figure}

This explanation aside, it is also of interest to numerically determine the so-called \textit{saturation point}: that is to say, the maximal value of $T$ (for a given $\kappa^*$) above which saturation occurs.  This is shown in Figure \ref{f:Saturation_Compute}, where we plot the quantity $\kappa^{(T)}_{M,M}$ against $T$ for $M=500$.  Note that this quantity is decreasing in $T$ (recall Figure \ref{f:ContourPlot}).  Using this figure, we make the following observation.  For the case $T=2$, saturation occurs when the $\kappa^* \approx 5e4$ or greater.  Hence, we now conclude the following: \textit{if one implements the equispaced FE with $T=2$, $\epsilon = 10^{-13}$ and $N = \Theta^{(2,\epsilon)}(M;\kappa^*)$ chosen such that the $\kappa^*$ is less than $\approx 5e4$, then no other value of $T$ will asymptotically give a better approximation.}  This establishes one of the main aims of this paper: namely, determining when the value $T = 2$, and hence the associated fast algorithm, can be used without concern that another value gives better accuracy.  Note that this conclusion is very reasonable.  In practice, we usually want the condition number to be much smaller than $10^4$ in magnitude.

Having now ascertained that one may use $T=2$ without worry in most cases, we end this subsection by providing numerical values for the constants $\kappa^{(2,\epsilon)}_{M / \eta , M}$ and $\lambda^{(2,\epsilon)}_{M/\eta,M}$.  These are shown in Figure \ref{f:T2_values}.  As we see, the choices $\eta = 1.5$ or $\eta=2.0$ seem reasonable in practice.  Increasing $\eta$ beyond this point brings only marginal benefits in stability.

\begin{figure}\centering
\begin{tabular}{|c|c|c|c|c|c|}
\hline $\eta/M$ & 250 & 500 & 750 & 1000  \\ \hline 
1.00  & 2.51e5 & 2.86e5 & 2.65e5 & 3.14e5   \\ \hline
1.25  & 1.01e4 & 1.25e4 & 1.72e4 & 1.99e4 \\ \hline
1.50  & 2.16e3 & 2.39e3 & 2.41e3 & 2.84e3 \\ \hline
2.00  & 1.88e2 & 2.25e2 & 2.89e2 & 3.27e2 \\ \hline
3.00  & 2.70e1 & 3.29e1 & 3.94e1 & 3.94e1 \\ \hline
4.00  & 1.18e1 & 1.53e1 & 1.67e1 & 1.84e1 \\ \hline
  \end{tabular}
    \hspace{1pc}
\begin{tabular}{|c|c|c|c|c|c|}
\hline $\eta/M$ & 250 & 500 & 750 & 1000  \\ \hline 
1.00  & 6.16e-7 & 1.48e-6 & 2.76e-6 & 3.49e-6  \\ \hline
1.25  & 3.86e-8 & 8.82e-8 & 9.98e-8 & 1.33e-7 \\ \hline
1.50  & 3.04e-9 & 7.36e-9 & 1.56e-8 & 1.82e-8 \\ \hline
2.00  & 4.26e-10 & 1.02e-9 & 1.00e-9 & 1.23e-9 \\ \hline
3.00  & 2.97e-11& 8.20e-11 & 8.37e-11 & 1.97e-10 \\ \hline
4.00  & 1.71e-11 & 1.51e-11 & 3.59e-11 & 5.24e-11 \\ \hline
  \end{tabular}
  \caption{Values of $\kappa^{(T,\epsilon)}_{M/\eta,M}$ (left) and $\lambda^{(T,\epsilon)}_{M//\eta,M}$ (right) for $T=2$ and $\epsilon = 10^{-13}$.} \label{f:T2_values}
\end{figure}

\subsection{Other solvers}\label{ss:other_solvers}

\begin{figure}
\begin{center}
$\begin{array}{ccc}
\includegraphics[width=5.00cm]{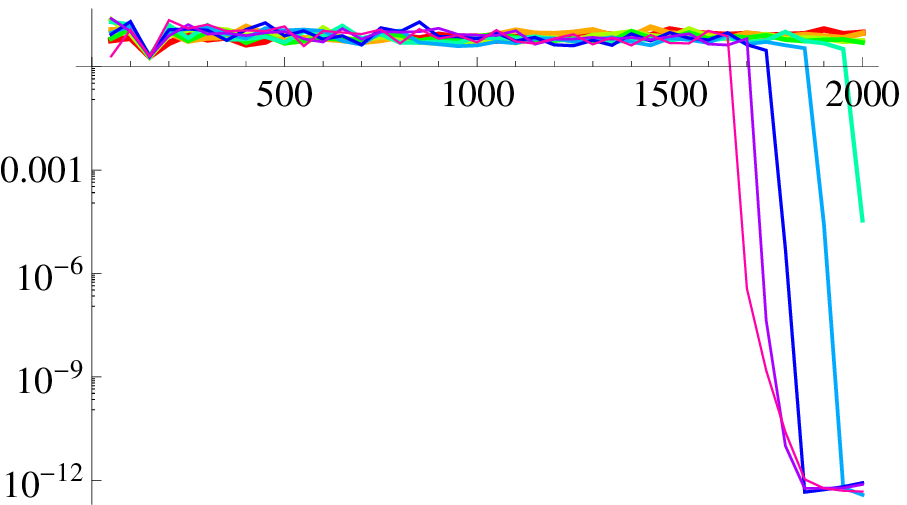} & \includegraphics[width=5.00cm]{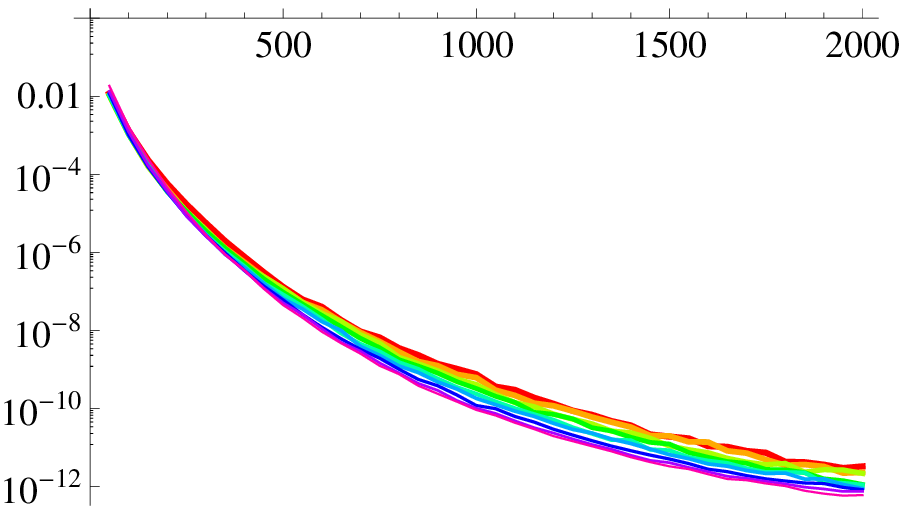} & \includegraphics[width=5.00cm]{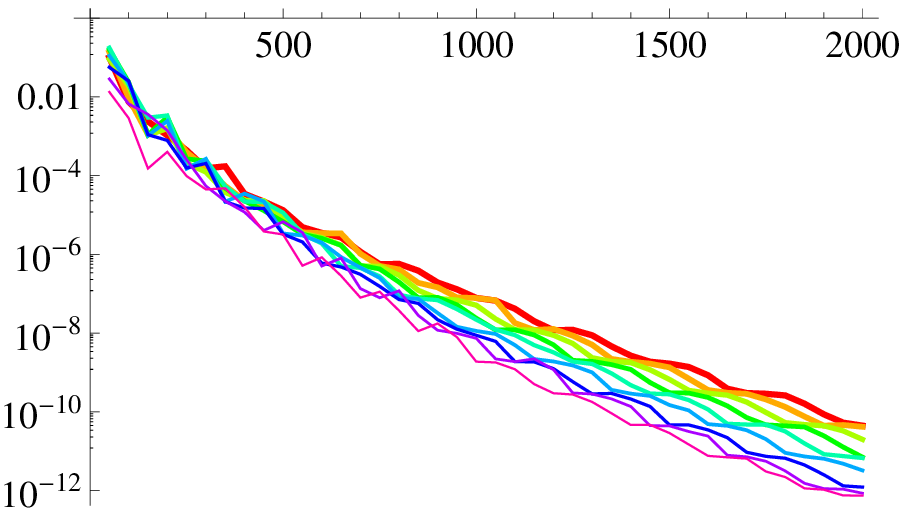} \\
\includegraphics[width=5.00cm]{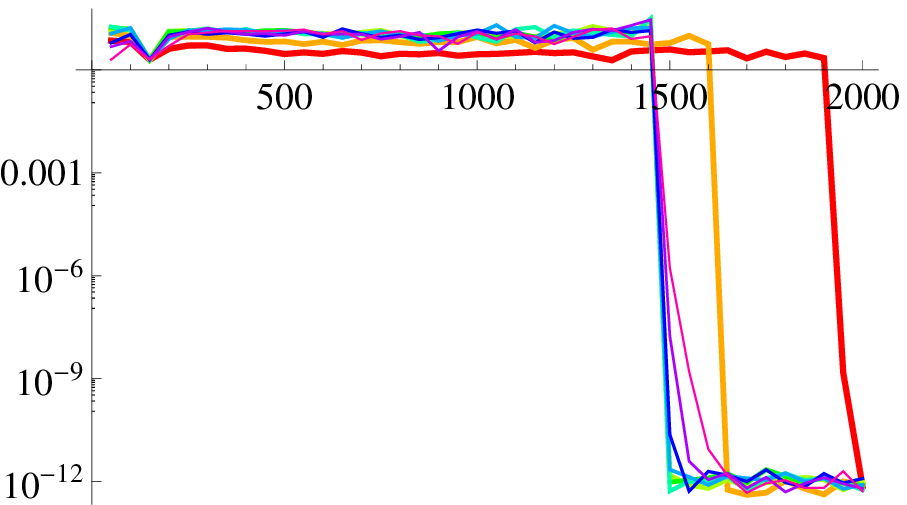} & \includegraphics[width=5.00cm]{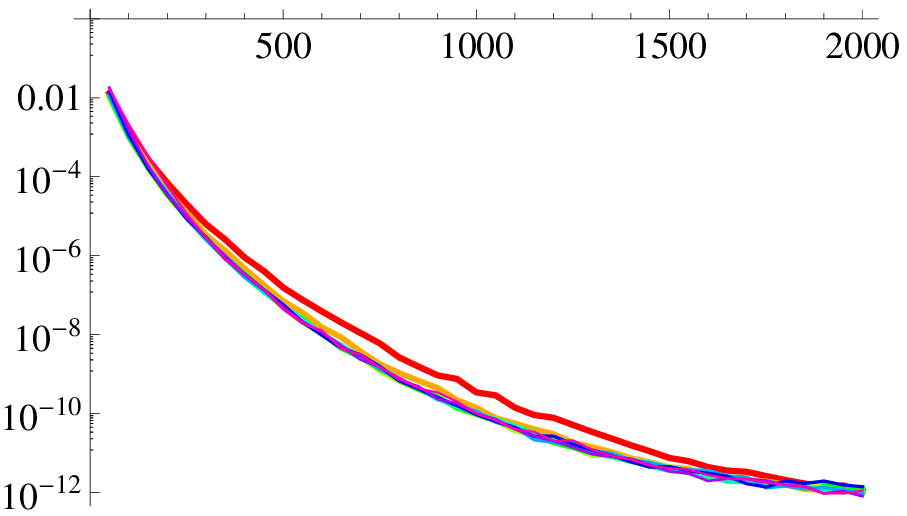} & \includegraphics[width=5.00cm]{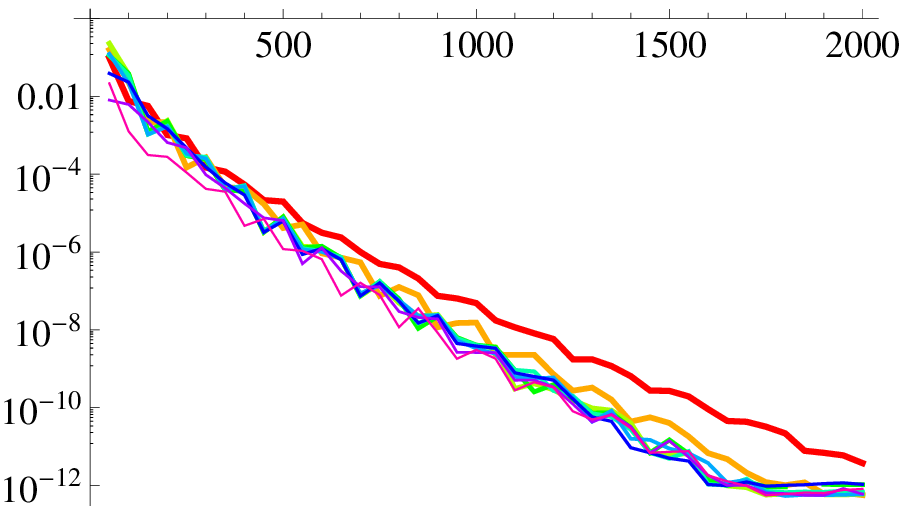}
\end{array}$
\caption{\small Approximation errors for the functions $f_1(x)$, $f_5(x)$ and $f_8(x)$ (left to right) using \textit{Matlab}'s $\backslash$ (top row) and \textit{Mathematica}'s \texttt{LeastSquares} (bottom row).  For each solver, the values $\Theta^{(T)}(M;\kappa^*)$ were computed for $\kappa^* = 25$.} \label{f:othersolvers}
\end{center}
\end{figure}

The phenomenon described above also occurs when a different numerical solver is used.  We illustrate this in Figure \ref{f:othersolvers} for \textit{Matlab}'s $\backslash$ and \textit{Mathematica}'s \texttt{LeastSquares}.  Interestingly, the effect is much more pronounced for the latter than for the former, although we believe it would become more apparent in the former for larger $M$ (recall that this phenomenon is asymptotic in $M$).  The discrepencies in the results are likely due to the different algorithms having somewhat different default tolerances for solving ill-conditioned least-squares problems.

\begin{figure}
\begin{center}
$\begin{array}{ccc}
\includegraphics[width=5.00cm]{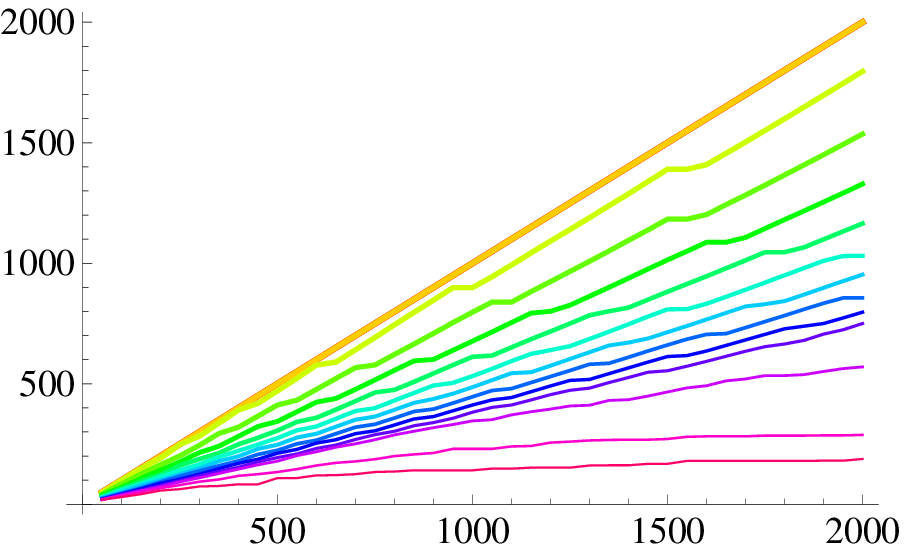}  & 
\includegraphics[width=5.00cm]{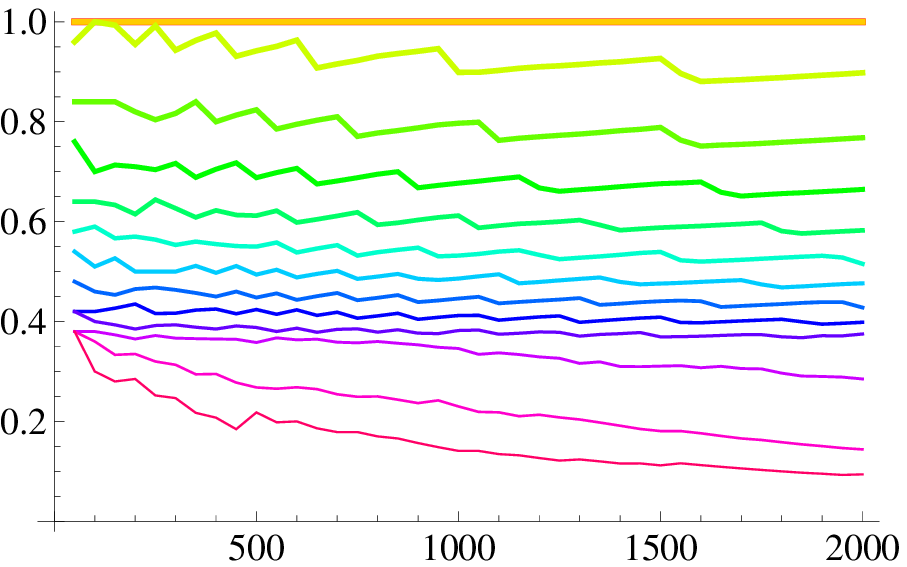}
&
\includegraphics[width=5.00cm]{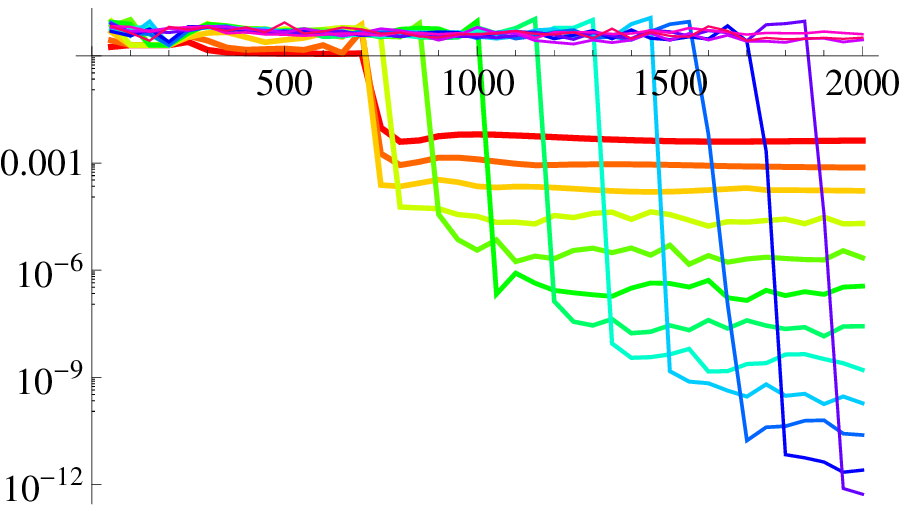}
\end{array}$
\caption{\small Left and Middle: the functions $\Theta^{(T,\epsilon)}(M;\kappa^*)$ and $\Theta^{(T,\epsilon)}(M;\kappa^*)/M$ against $M$ for $T=2$, $\kappa^* = 10$ and $\epsilon = 10^{-2},10^{-3},\ldots,10^{-16}$ (thickest to thinnest).  Note that the graphs for $\epsilon = 10^{-2},10^{-3},10^{-4}$ are the same.  Right: approximation of $f(x) = \exp(250 \sqrt{2} \pi \I x)$ using these values.} \label{f:FnAppLog}
\end{center}
\end{figure}

\subsection{Influence of the SVD tolerance $\epsilon$}

Let us now return to the SVD algorithm.  Thus far, we have taken the tolerance $\epsilon$ to be equal to $10^{-13}$.  However, other values are possible, and it is of interest to determine how this affects the approximation.  Note that some previous insight in this problem was given in \cite{BoydFourCont}, where the effect of $\epsilon$ was considered for specific functions.

In Figure \ref{f:FnAppLog} we plot $\Theta^{(T,\epsilon)}(M;\kappa^*)$ for different values of $\epsilon$ using $T=2.0$.  As is evident, $\Theta^{(T,\epsilon)}(M;\kappa^*)$ grows more rapidly for increasing $\epsilon$.  This should come as no surprise, since it is the small singular vectors that cause ill-conditioning.  However, the main purpose of this figure is to show the specific improvement that is possible by changing $\epsilon$.  For example, when $\epsilon = 10^{-13}$ the function $\Theta^{(T,\epsilon)}(M;\kappa^*)$ is approximately $0.4 M$.  Conversely, for $\epsilon = 10^{-6}$ it is approximately $0.8 M$, i.e.\ it scales roughly twice as quickly.  If $\delta$ is some finite tolerance greater than $10^{-6}$, this means that the FE approximation with $\epsilon = 10^{-6}$ will approximate a given function to accuracy $\delta$ using roughly half the number of equispaced points as is required when $\epsilon = 10^{-13}$.  This fact is also confirmed in the right panel of Figure \ref{f:FnAppLog}, where the oscillatory function $f(x) = \exp(250 \sqrt{2} \pi \I x)$ is approximated using different values of $\epsilon$.  Resolving this function using $\epsilon = 10^{-13}$ requires roughly $3700$ equispaced data points, whereas when $\epsilon = 10^{-6}$ this value drops to around $1700$.

Of course, the downside of a larger $\epsilon$ is that the minimal error is limited to approximately $\epsilon$, as can be seen in Figure \ref{f:FnAppLog}. Nonetheless, the conclusion we draw from this section is that if accuracy close to machine precision is not required -- as is typically the case in practice, where three to six digits is often acceptable -- then a viable way to increase the performance of the FE approximation whilst maintaining the condition number is to use a larger value of $\epsilon$.

\section{Resolution power}

In many applications, it is important to have an approximation algorithm with good \textit{resolution power}.  Loosely speaking, this means that oscillatory functions are recovered using using a number of measurements that scales linearly with the frequency of oscillation with a constant that is as small as possible.  Formally, let $\{ F_M \}_{M \in \bbN}$ be a sequence of approximations such that $F_M(f)$ depends only on the values of $f$ on an equispaced grid of $2M+1$ points.  Let
\bes{
\cR(\omega,\delta) = \min \left \{ M \in \bbN : \| \E^{ \I  \pi \omega \cdot} - F_M(\E^{\I \pi \omega \cdot}) \|_{\infty} < \delta \right \},\quad \omega > 0,\  0 < \delta < 1,
}
then we say that $F_M$ has resolution constant $0 < r <\infty$ if 
\be{
\label{}
\cR (\omega,\delta) \sim r \omega,\quad \omega \rightarrow \infty,
}
for any fixed $\delta$.  The approximation $F_M$ has good resolution power if $r$ is small, and bad resolution power if $r$ is large.

For periodic oscillations, the Fourier series approximation on $[-1,1]$ of degree $M$ has optimal resolution constant $r=1$.  Since the number equispaced points is $2M+1$, this corresponds to two \textit{points-per-wavelength}.  Of course, Fourier series do not converge uniformly for nonperiodic oscillations, which is why we resort to alternative algorithms such as FEs.  Naturally, though, it is desirable that the FE resolution constant be as close to the optimal value $r=1$ as possible.

\begin{figure}
\begin{center}
$\begin{array}{ccc}
\includegraphics[width=5.00cm]{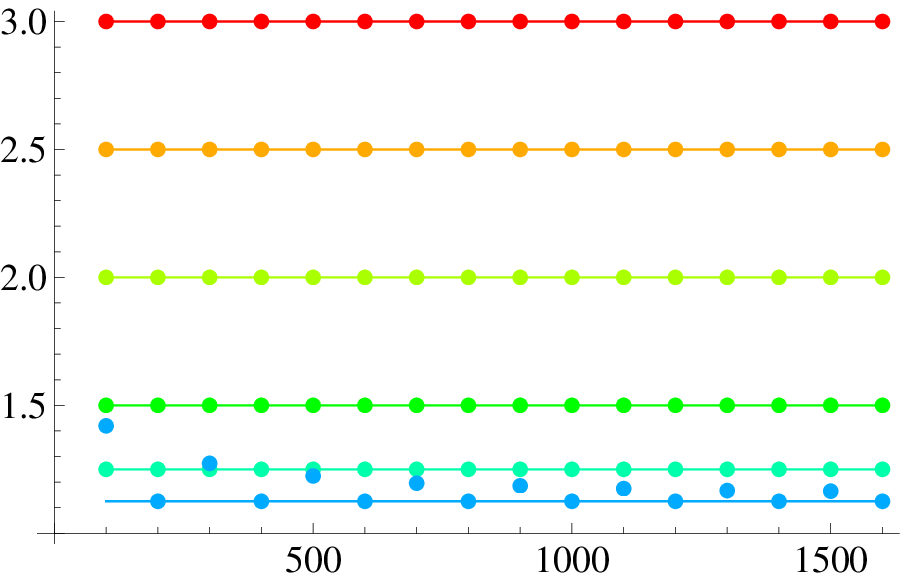} & 
\includegraphics[width=5.00cm]{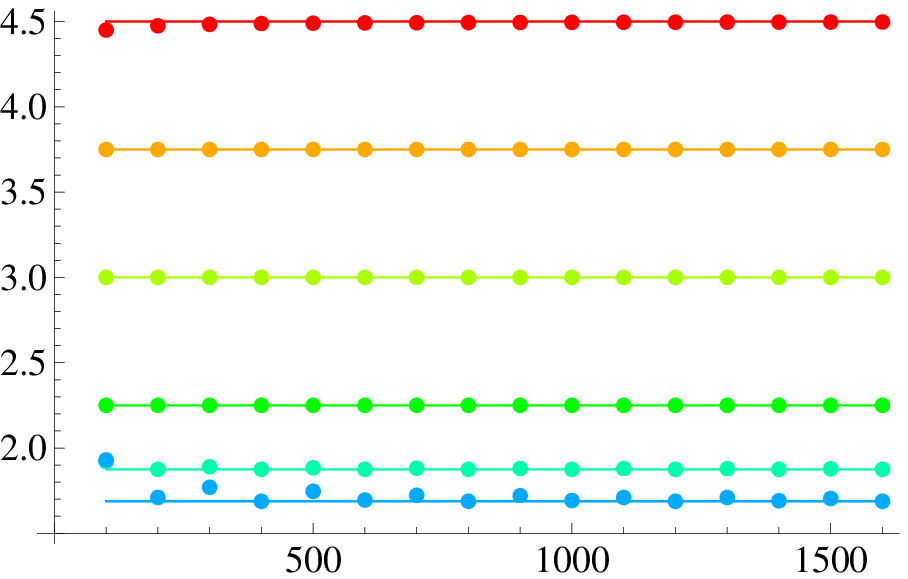} & 
\includegraphics[width=5.20cm]{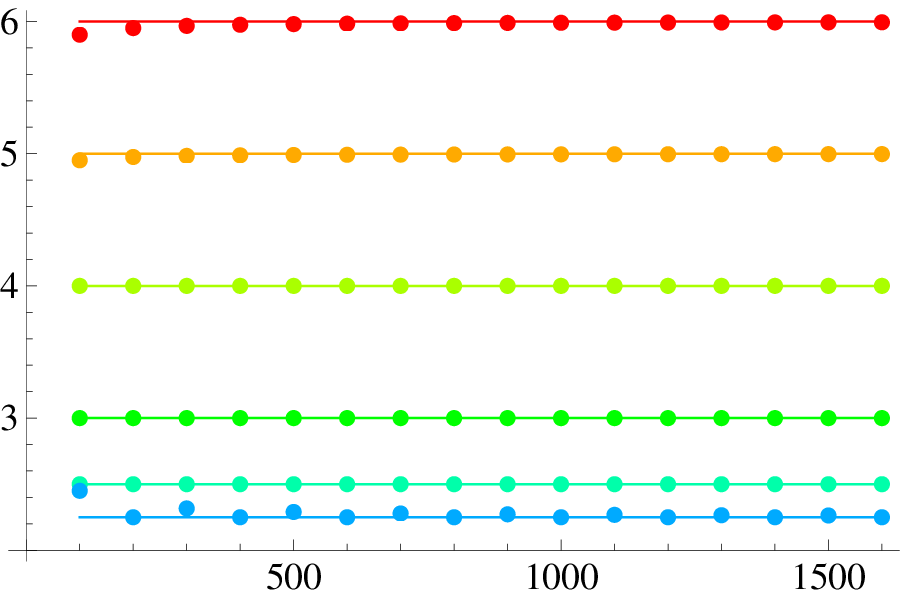} \\
\eta = 1 & \eta = 1.5 & \eta = 2
\end{array}$
\caption{\small The quantity $\cR(\omega,\delta ; T,\eta) / \omega$ against $\omega$ for $T=3.0,2.5,2.0,1.5,1.25,1.125$ and $\delta = 10^{-3}$.  The solid lines indicate the values $T \eta$.} \label{f:ResnPower}
\end{center}
\end{figure}

In Figure \ref{f:ResnPower} we numerically determine the resolution constant for the FE approximation $F^{(T,\epsilon)}_{M/\eta,M}$ by computing the function
\bes{
\cR(\omega,\delta)  =\cR(\omega,\delta ; T , \eta).
}
These results suggest that the resolution constant
\be{
\label{r_prod}
r \approx T \eta,
}
is approximately the product of the extension parameter $T$ and the oversampling ratio $\eta$.

We shall discuss the consequences of this observation in a moment, but we first wish to explain why $r$ should be at most $T \eta$.  Recall from \R{error_numerical_2} that the error of the FE approximation is determined (up to a small linear constant in $M$) by the decay of the factor $E_N(f)$, given by \R{EN_def}, where $N = M/\eta$.  Now let
\be{
\label{R0}
\cR_0(\omega,\delta) = \min \left \{ N \in \bbN : E_N(\E^{\I \pi \omega \cdot})  < \delta \right \},\quad \omega > 0, \sqrt{2T} \mu < \delta < 1.
}
Then to show that $r = T \eta$, it suffices to prove that
\bes{
\cR_0(\omega,\delta) \sim T \omega,\quad \omega \rightarrow \infty.
}
We now note the following:

\lem{
Let $\cR_0(\omega,\delta)$ be as in \R{R0}.  Then
\bes{
\limsup_{\omega \rightarrow \infty} \cR_0(\omega,\delta)/\omega \leq T.
}
}
\prf{
Let $f(x) = \E^{\I \pi \omega x}$ and write $\omega = P/T + z$ where $0 \leq z < 1/T$ and $P \in \bbN$.  Then
\bes{
f(x) = \E^{\I \frac{P \pi}{T} x} g(x),\qquad g(x) = \E^{\I \frac{z \pi}{T} x}.
}
Let $K \in \bbN$ and let $\phi \in \cG^{(T)}_{P+K}$ be given by $\phi(x) = \E^{\I \frac{P \pi}{T} x} \tilde{\phi}(x)$, where $\tilde{\phi} \in \cG^{(T)}_{K}$ is arbitrary.  Define $N = P+K$, and let $\mathbf{a}$ be the vector of coefficients of $\phi$.  Then by definition,
\bes{
E_N(f) \leq \| f -\phi \|_{\infty} + \mu | \mathbf{a} |_{\infty} = \| g - \tilde{\phi} \|_{\infty} + \mu | \mathbf{\tilde{a}} |_{\infty},
}
where $\mathbf{\tilde{a}}$ is the vector of coefficients of $\tilde{\phi}$.  Since $\tilde{\phi}$ is arbitrary, we deduce that $E_N(f) \leq E_{K}(g)$.  Thus $E_N(f) \leq \delta$ provided $E_K(g) \leq \delta$.  Using Corollary \ref{c:EN_rate}, we see that $E_K(g) \leq \delta$ provided $K \geq K_0(\delta)$, where $K_0(\delta)$ is independent of $\omega$.  Thus, $E_N(f) \leq \delta$ whenever $N = P + K_0(\delta) = (\omega - z)T + K_0(\delta)$.  Since $K_0(\delta)$ is independent of $\omega$, we obtain the result.
}
From this, we deduce that $\cR_0(\omega,\delta) \leq T \omega + \ord{1}$ for large $\omega$ and hence $r \leq T \omega$.  Unfortunately, we have no proof of the lower bound, although it is supported by the results in Figure \ref{f:ResnPower}.

Let us now discuss the consequences of \R{r_prod}.  Suppose we return to the earlier experiment where $\kappa^*$ is fixed and $N$ is chosen according to $N = \Theta^{(T)}(M;\kappa^*)$.  For large $M$, Figure \ref{f:ThetaPlot} shows that
\be{
\label{etaTkappa}
\Theta^{(T)}(M;\kappa^*) \approx M/\eta^{(T)} ,\qquad \eta^{(T)} = \eta^{(T)}(\kappa^*),
}
for some fixed $\eta^{(T)}$ depending on $T$ and $\kappa^*$.  In our previous observation, it was found that the resulting FE approximation was independent of $T$.  In particular, this holds for the oscillatory exponential $\exp(\I \pi \omega x)$.  But since the resolution constant $r$, which describes the point after which the approximation error begins to decay, is equal to $\eta T$, we therefore deduce the following:
\bes{
\eta^{(T)} T = \eta^{(T')} T',
}
for any $T$ and $T'$ that do not saturate.  This implies that the level curves of the condition number $\kappa = \kappa^{(T)}_{M/\eta,M}$ are approximately given by $T \eta = \mbox{constant}$, for sufficiently large $M$ and non-saturating $T$, which is in good agreement with the contour plot given in Figure \ref{f:ContourPlot}.

In addition to this, another important implication of \R{r_prod} is that to get better resolution power one must necessarily worsen the stability of the algorithm.  In other words, there is a direct relationship between $\kappa^*$ and $r$, regardless of the choice of $T$ made.  In Table \ref{tab:ResnConstT2} we give numerical results for the resolution constant when $T=2$ (and therefore all nonsaturating $T$) for different values of $\kappa^*$.  As we see, by allowing $\kappa^*$ to increase to roughly $500$ we get a marked improvement over when $\kappa^* = 10$.  Beyond this point, further increases give only marginal gains at the expense of much larger condition numbers.

\begin{table}
\begin{center}
\begin{tabular}{|c|c|c|c|c|c|c|c|c|c|}
\hline
$\kappa^*$ & $10$ & $25$ & $100$ & $500$ & $1000$ & $5000$ & $10000$  \\
\hline
$r$ & 5.41 & 4.44 & 3.64 & 2.98 & 2.77 & 2.38 & 2.25\\
\hline
\end{tabular}
\caption{The resolution constant $r = \eta^{(T)} T$ for $T=2$, where $\eta^{(T)} = \eta^{(T)}(\kappa^*)$ is given by \R{etaTkappa}.}\label{tab:ResnConstT2}
\end{center}
\end{table}

\section{Other data}\label{s:other_data}
In this final section, we illustrate that the main phenomena observed for equispaced data are also witnessed for numerous other types of data.  To this end, we consider the following four examples:

\begin{itemize}
\item \textit{Jittered pointwise data}.  Here the measurements of $f$ are pointwise samples at the jittered locations
\bes{
f(x_m),\qquad x_m = \frac{m}{M} + z_m,
}
where $z_m \in (-\delta/M,\delta/M)$ and $0 < \delta < 1$.  This is a typical example of a nonuniform sampling pattern for scattered data approximation.  We shall choose the $z_m$'s as follows:
\be{
\label{jittered_data}
z_m = \frac{\delta}{M} \sin(M^2/m),\quad m \neq 0,\quad z_0 = 0.
}
\end{itemize}

\begin{itemize}
\item \textit{Logarithmic pointwise data}.  Here we again sample $f$ pointwise at nodes $x_m$, but in this case the nodes are logarithmically distributed:
\be{
\label{log_data}
x_m = - x_{-m} = 10^{\left(\frac{m-1}{M-1} - 1\right) \log_{10}(c M)},\quad x_0 = 0,
}
where $c>0$ is a fixed, user-controlled parameter (we take $c=2$ in our results).  This sampling pattern corresponds to a nonuniform sampling scenario where data is collected more densely at the origin.
\end{itemize}

\begin{itemize}
\item \textit{Fourier data}. In some applications, rather than pointwise samples, we may wish to reconstruct $f$ from its Fourier coefficients:
\be{
\label{Four_data}
\hat{f}(m) = \int^{1}_{-1} f(x) \E^{-\I \pi m x} \D x,\quad |m| \leq M.
}
As discussed in \cite{BAACHOptimality}, this can be seen as a continuous analogue of the equispaced data recovery problem.  In particular, there is a completely analogous result result to that of Platte, Trefethen \& Kuijlaars regarding stability and convergence \cite{AdcockHansenShadrinStabilityFourier}.  Note that for this data we replace the uniform norms used in the error estimates and the definitions of $\kappa^{(T)}_{N,M}$ and $\lambda^{(T)}_{N,M}$ by the $\rL^2$ and $\ell^2$ norms.  This is natural in view of Parseval's identity for Fourier coefficients.
\end{itemize}

\begin{itemize}
\item \textit{Optimal pointwise data}. Finally, in order to show that the phenomenon is not witnessed for all data, we consider pointwise samples taken at the so-called mapped symmetric Chebyshev nodes (see \cite{BADHFEResolution}):
\be{
\label{mapped_Cheby}
x_m = - x_{-m-1} = m^{-1} \left ( \cos \frac{(2m+1) \pi}{2M+2} \right ),\ m=0,\ldots,M,
}
where $m$ is given by \R{mapping}.  These nodes are derived from the observation that FE approximations correspond to algebraic polynomial approximations in the mapped co-ordinate $z = m(x)$ \cite{FEStability}.  Chebyshev nodes provide optimal nodes for polynomial interpolation.  Therefore, under the inverse mapping $m^{-1}$, they provide the optimal nodes \R{mapped_Cheby} for FE approximations.  Note that since the nodes arise in this way, no oversampling is required in the FE approximation, i.e.\ we let $\eta = 1$ in this case.
\end{itemize}

\begin{figure}
\begin{center}
$\begin{array}{ccc}
\includegraphics[width=5.00cm]{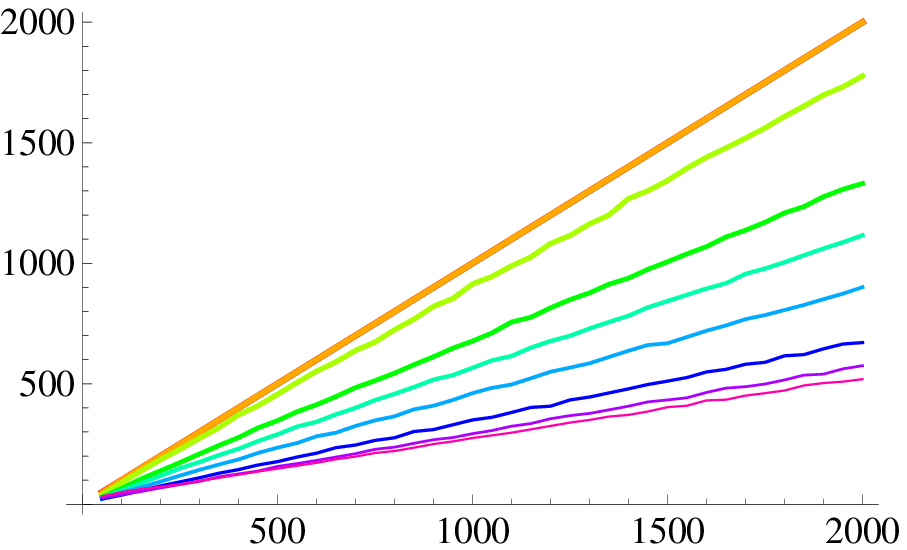}  &  \includegraphics[width=5.00cm]{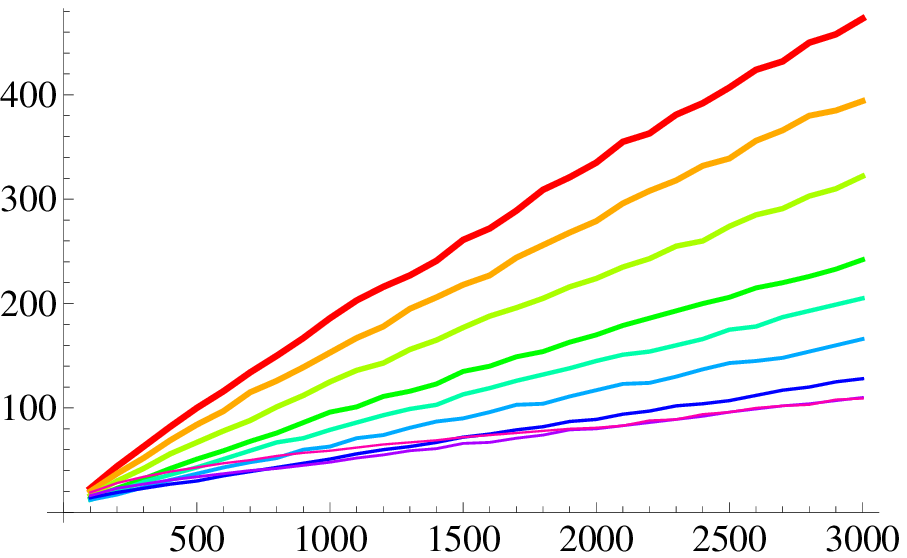} & \includegraphics[width=5.00cm]{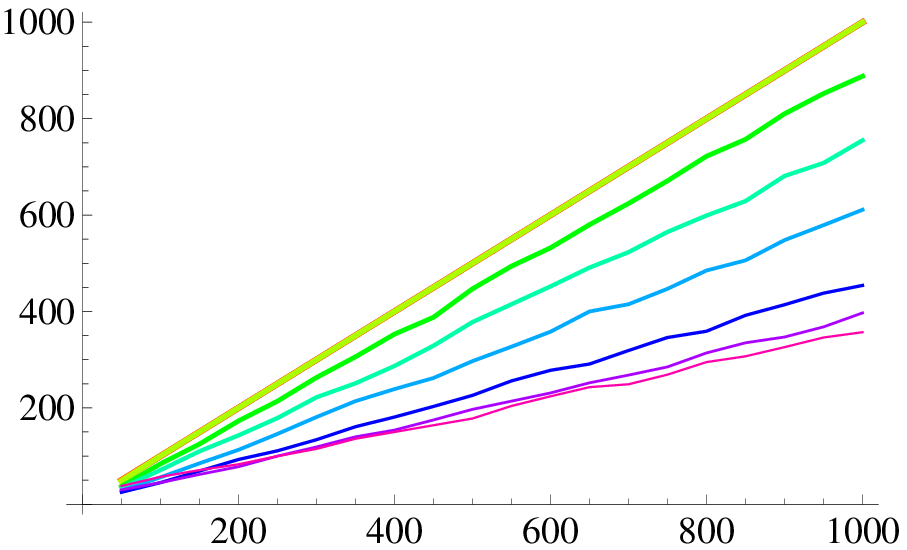}
\\
\includegraphics[width=5.00cm]{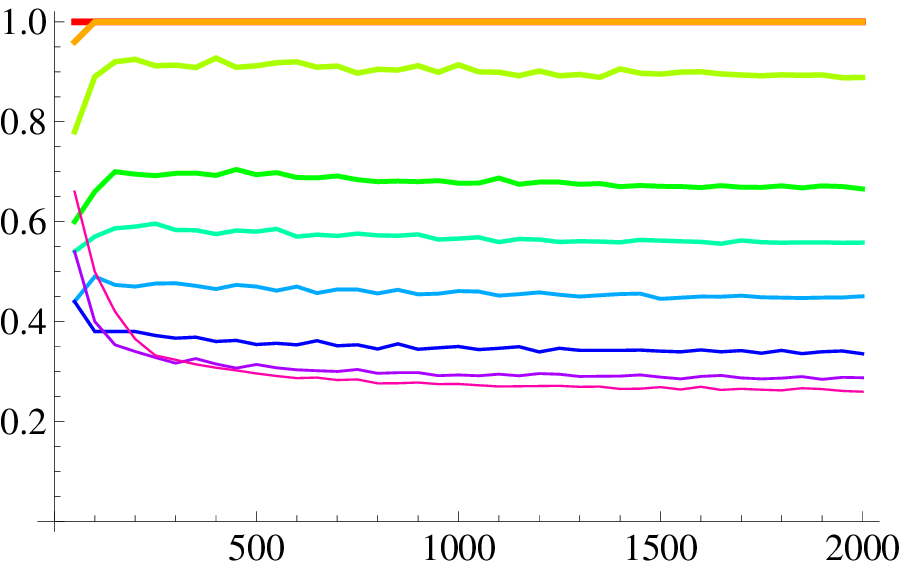} & \includegraphics[width=5.00cm]{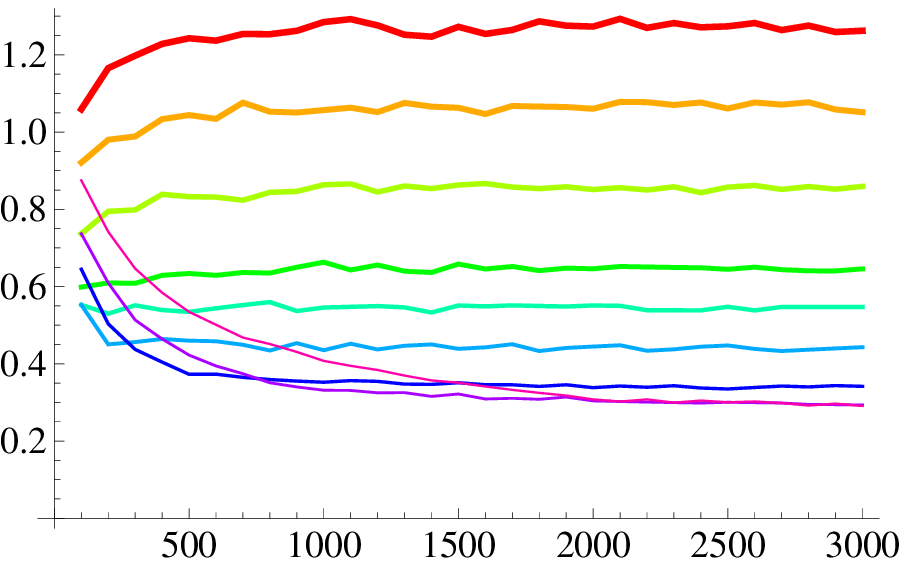} 
& \includegraphics[width=5.00cm]{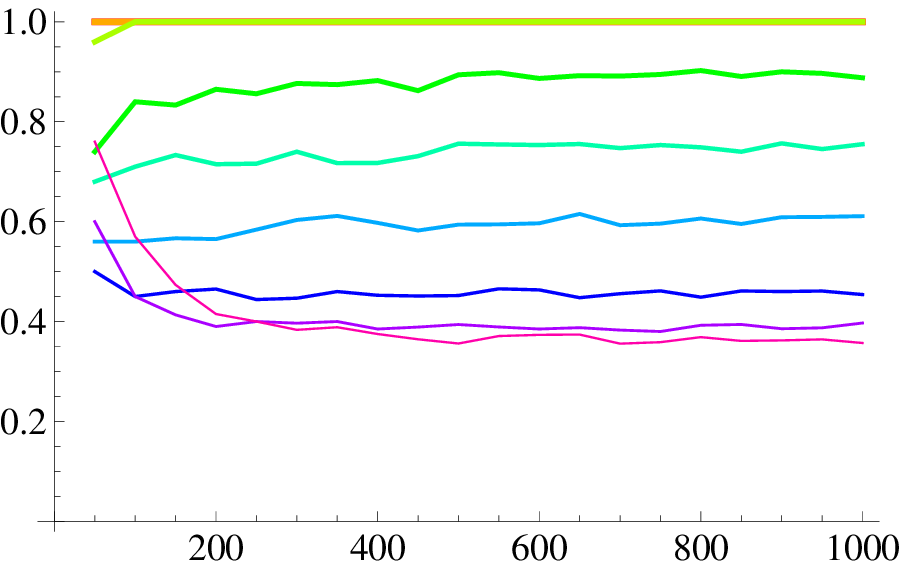}
\end{array}$
\caption{\small Plots of $\Theta^{(T)}(M;\kappa^*)$ (top row) and $\Theta^{(T)}(M;\kappa^*) / S(M)$  (bottom row) against $M$ for jittered \R{jittered_data} (left), logarithmic \R{log_data} (middle) and Fourier \R{Four_data} (right) data.  Here $S(M) = M/\log(M)$ for the logarithmic data and $S(M) = M$ otherwise, and $\kappa^* = 25$ (jittered, logarithmic) or $\kappa^* = 10$ (Fourier).} \label{f:ThetaPlot_Data}
\end{center}
\end{figure}

In Figure \ref{f:ThetaPlot_Data} we give plots of the function $\Theta^{(T)}(M;\kappa^*)$ for the first three data types.  For jittered and Fourier data the scaling is linear, whereas for the logarithmic data $\Theta^{(T)}$ scales like $M / \log M$.  This scaling is proportional to the reciprocal of the maximal spacing between nodes in the case, and hence is completely expected.  Note also that no values of $T$ saturate for the logarithmic data, whereas values $T=5.0$ and $T=6.0$ saturate for the jittered data, and  for the Fourier data the values $T=3.0$, $T=4.0$, $T=5.0$ and $T=6.0$ all saturate.  The lower saturation point for the latter is due to the fact that the condition number is measured in the weaker $\rL^2$ norm.

\begin{figure}
\begin{center}
$\begin{array}{ccc}
\includegraphics[width=5.20cm]{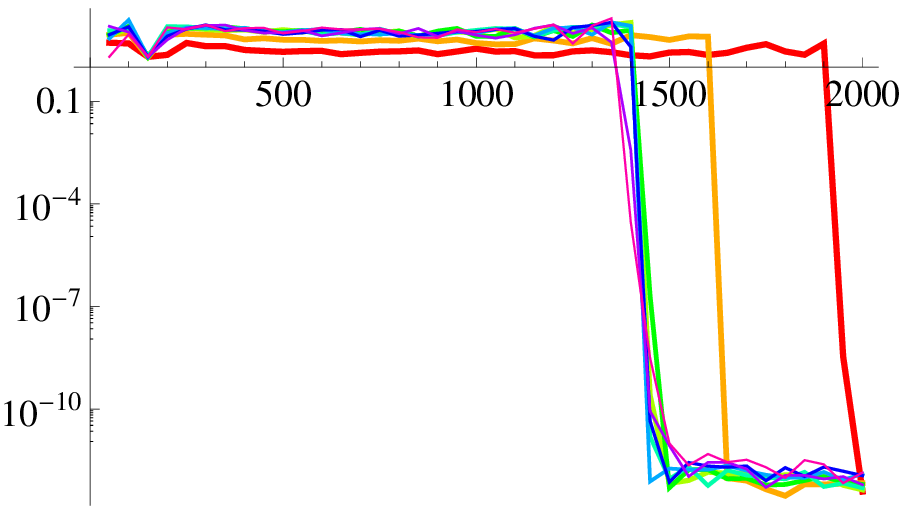}  
& 
\includegraphics[width=5.20cm]{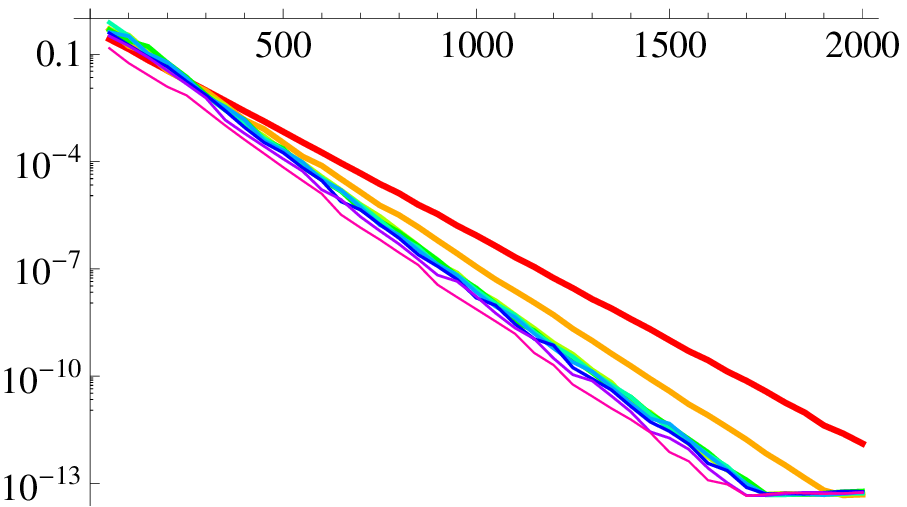}  
& 
\includegraphics[width=5.20cm]{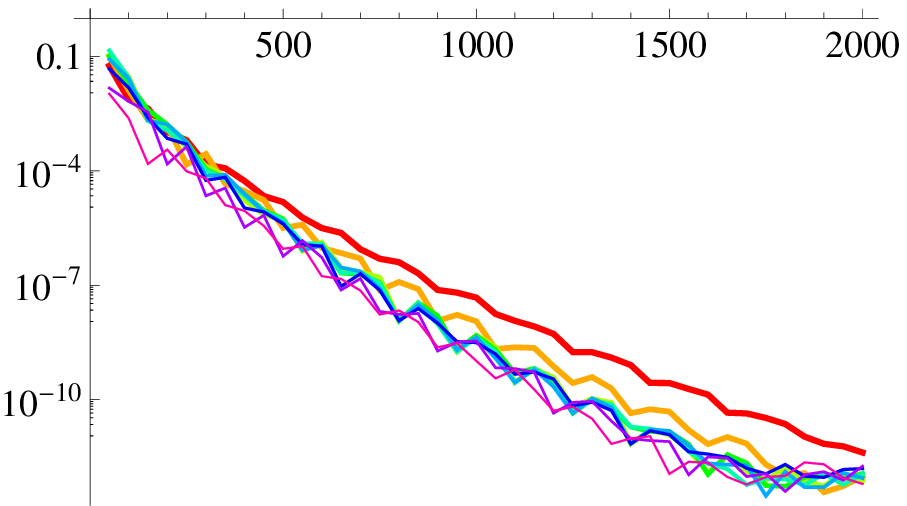}  
\\
f(x) = \E^{230 \sqrt{2} \pi \I x} & f(x) = \frac{1}{1+1500 x^2} & f(x) = \E^{-1/(8 x)^2}
\end{array}$
\caption{\small Approximation errors for jittered data \R{jittered_data} using the values $\Theta^{(T)}(M;\kappa^*)$ from Figure \ref{f:ThetaPlot_Data}.}  \label{f:FnApp_Data}
\end{center}
\end{figure}

\begin{figure}
\begin{center}
$\begin{array}{ccc}
\includegraphics[width=5.20cm]{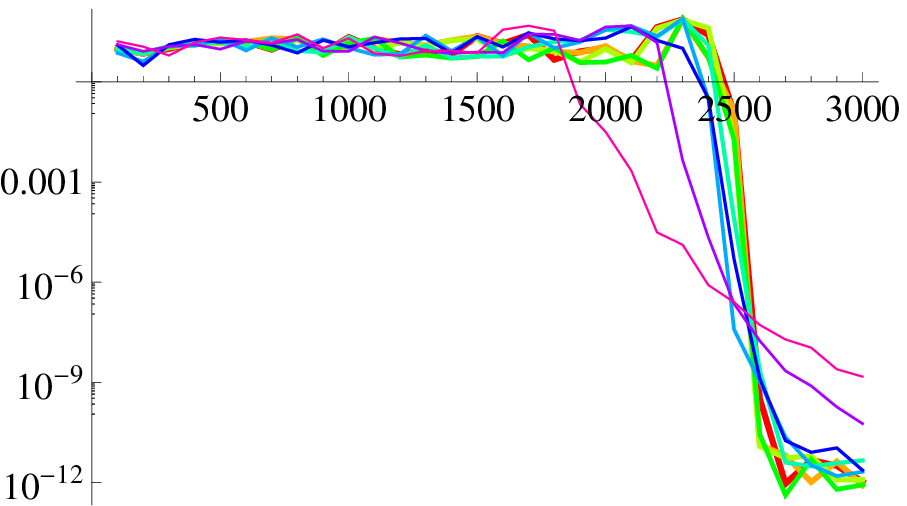}  
& 
\includegraphics[width=5.20cm]{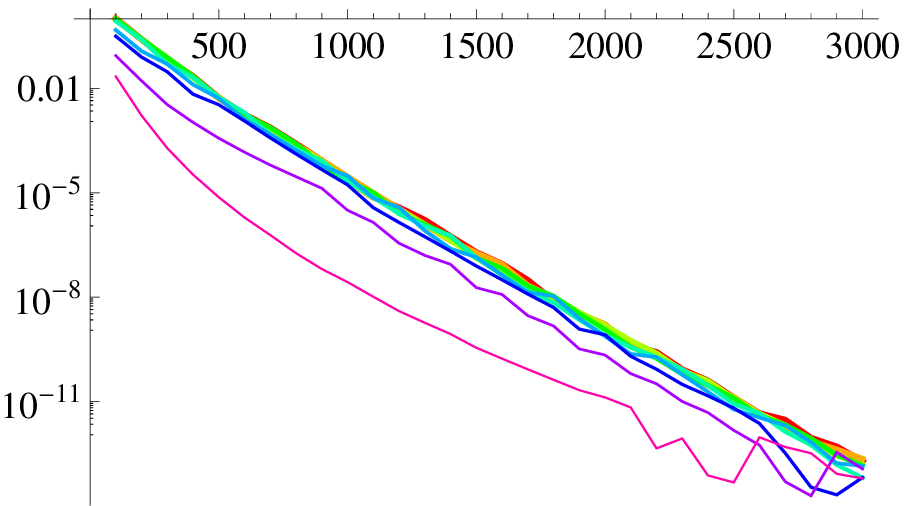}  
& 
\includegraphics[width=5.20cm]{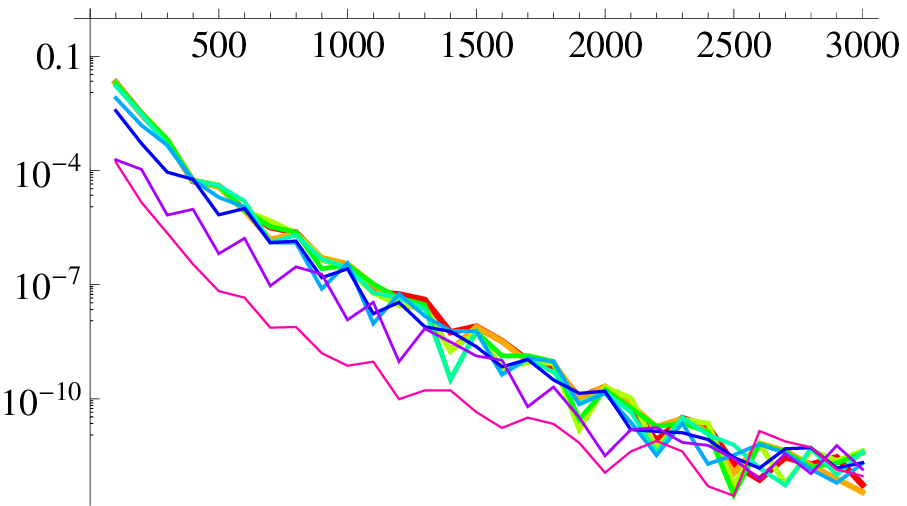}  
\\
f(x) = \E^{50 \sqrt{2} \pi \I x} & f(x) = \frac{1}{1+65 x^2} & f(x) = \E^{-2/(3x^2)}
\end{array}$
\caption{\small Approximation errors for logarithmic data \R{log_data} using the values $\Theta^{(T)}(M;\kappa^*)$ from Figure \ref{f:ThetaPlot_Data}. } \label{f:FnApp_Data2}
\end{center}
\end{figure}

\begin{figure} 
\begin{center}
$\begin{array}{ccc}
\includegraphics[width=5.00cm]{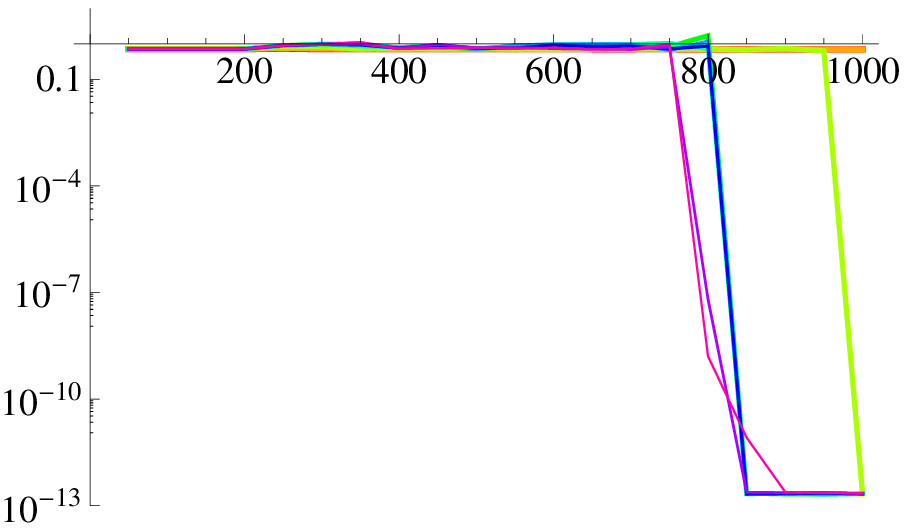}  &
\includegraphics[width=5.00cm]{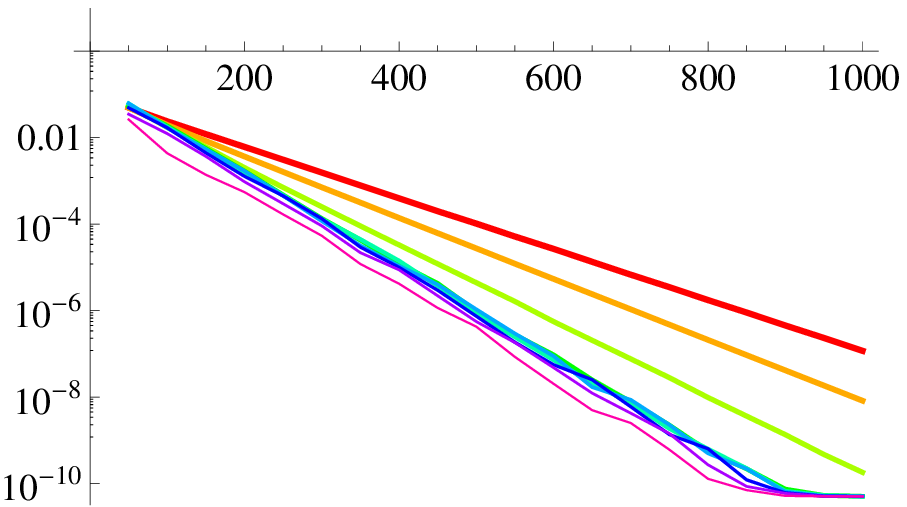}  &
\includegraphics[width=5.00cm]{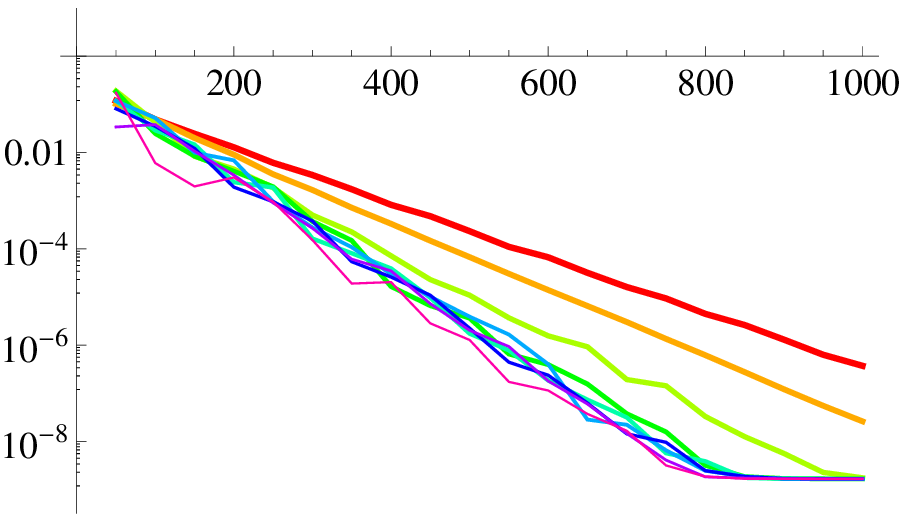}  
\\
f(x) = \E^{175 \sqrt{2} \pi \I x} & f(x) = \frac{1}{1+1500 x^2} & f(x) = \frac{1}{1+25 \sin^2 8x}
\end{array}$
\caption{\small Approximation errors for Fourier data \R{Four_data} using the values $\Theta^{(T)}(M;\kappa^*)$ from Figure \ref{f:ThetaPlot_Data}. } \label{f:FnApp_Data3}
\end{center}
\end{figure}

\begin{figure}
\begin{center}
$\begin{array}{ccc}
\includegraphics[width=5.00cm]{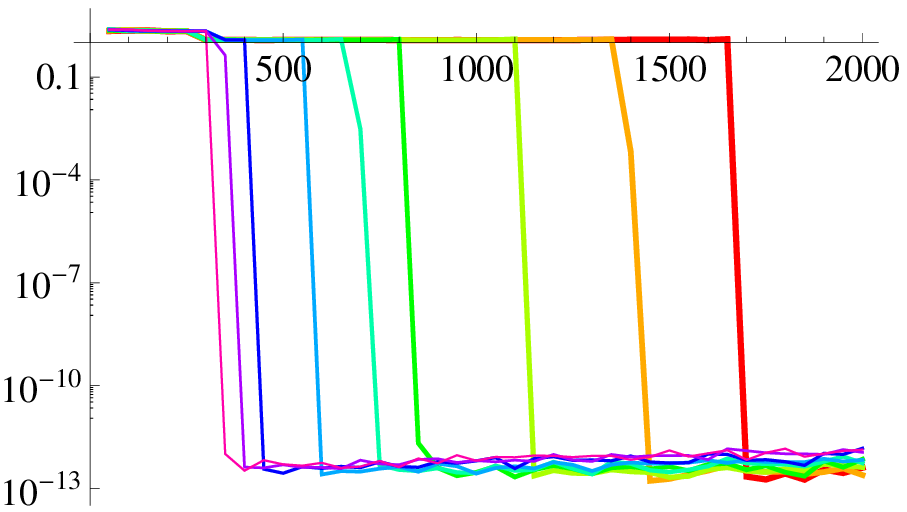}  & \includegraphics[width=5.00cm]{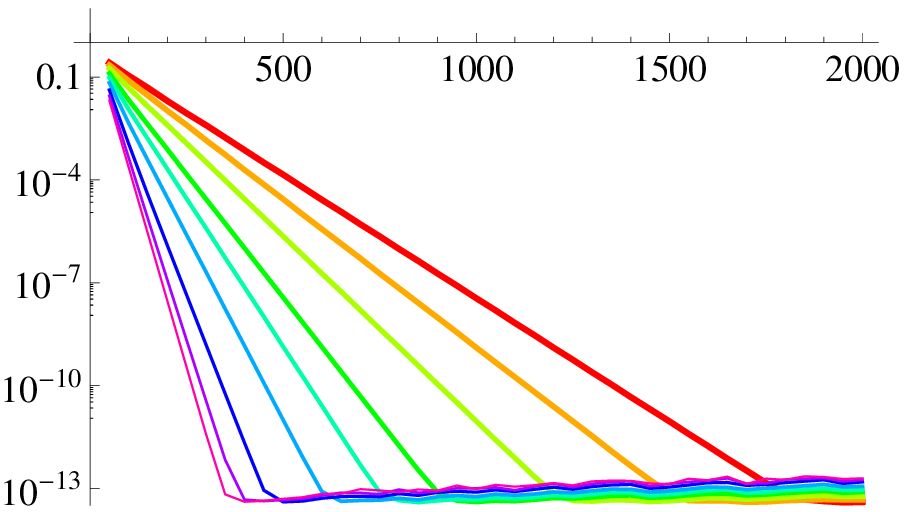} &
\includegraphics[width=5.00cm]{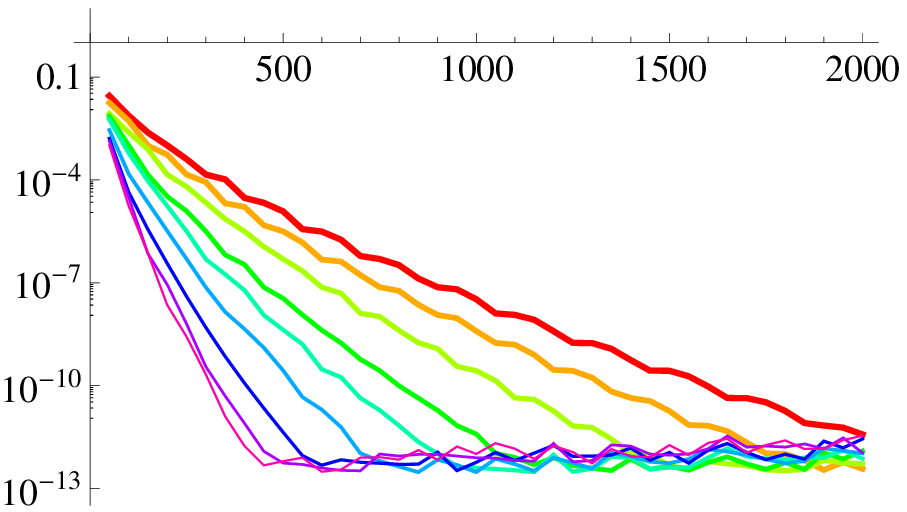} \\
f(x) = \E^{200 \sqrt{2} \pi \I x} & f(x) = \frac{1}{1+1000 x^2} & f(x) = \E^{-1/(8x)^2}
\end{array}$
\caption{\small Approximation errors for the optimal pointwise data \R{mapped_Cheby} using $M=N$.} \label{f:FnApp_Data4}
\end{center}
\end{figure}

Next, in Figures \ref{f:FnApp_Data}--\ref{f:FnApp_Data3} we compare approximation errors using these values.  For the jittered and Fourier data we see exactly the same phenomenon as before: namely, the approximation errors are roughly independent of $T$.  A similar phenomenon is witnessed for the logarithmic data, although it is slightly weaker: $T=1.125$ gives somewhat better errors than the other values.  This is due to the much more severe scaling of $\Theta^{(T)}(M;\kappa^*)$ with $M$ in this case, which means the asymptotic regime takes longer to set in.

Finally, in Figure \ref{f:FnApp_Data4} we display approximation errors for the optimal pointwise data \R{mapped_Cheby}.  As is evident, the phenomenon does not occur in this case.  The reason for this is due to the choice of the data, which means that no oversampling is required.  Thus the FE approximation error is proportional to $E_N(f)$, where $N = M$, multiplied by a mildly growing factor.  Figure \ref{f:FnApp_Data4} therefore serves as a reminder that smaller values of $T$ are intrinsically better than larger values, provided one has freedom to pick ideal data.  On the other hand, for nonideal data -- such as equispaced, jittered, logarithmic or Fourier data -- this effect is nullified by a worse scaling of $\Theta^{(T)}(M;\kappa^*)$.

\section{Conclusions and open problems}\label{s:conclusion}
The purpose of this paper was to document an interesting phenomenon in FE approximations from equispaced data.  Namely, when the desired condition number is fixed, the choice of the extension parameter $T$ has no substantial effect on the approximation.  This is on the proviso that saturation does not occur, which we have shown to be the case for moderate values of $T$ and $\kappa^*$.  In particular, one may use $T=2$, and the associated fast algorithm, without concern that it is suboptimal.

The main open problem is to provide mathematical analysis for the empirical conclusions drawn.  We believe this is possible, although not straightforward.  One possible approach towards this is to conduct an asymptotic analysis of the singular values and vectors of the matrix $A^{(T)}$.  Recall that the normal form $(A^{(T)})^* A^{(T)}$ is a Riemann sum approximation to the Gram matrix $G$ of the FE basis functions $\phi_n(x) = \E^{\I \frac{n \pi}{T} x}$.  As discussed in \cite{FEStability}, the matrix $G$ is precisely the prolate matrix.  The eigenvalues and eigenvectors of this matrix were analyzed in detail by Slepian \cite{SlepianV} (see also \cite{Varah}).  It may be possible to do the same for the discretized version  $(A^{(T)})^* A^{(T)}$, and this is an important topic for future work.

Another question raised by this work is that of whether it might be possible to vary $T$ with $M$ to achieve better results; in particular, improved resolution power.  We believe this may be the case, the caveat being that there is  currently no fast algorithm for $T \neq 2$.  Some potential choices for varying $T$ with $M$ were considered previously in \cite{BADHFEResolution,FEStability}.  But it may also be possible using the approach of this paper to numerically compute an optimal (in some sense) value of $T$ for each $M$.  We leave this for future work.

\section*{Acknowledgments}
The authors would like to thank Daan Huybrechs, Mark Lyon and Rodrigo Platte for useful discussions. BA acknowledges support from the NSF DMS grant 1318894.

\bibliographystyle{abbrv}
\small
\bibliography{FEParameterBib}

\newpage
\section*{Symbols}

\small

\renewcommand{\arraystretch}{1.2}
\begin{longtable}{|c|p{9.0cm}|}
 \hline Symbol &  
Description
\\ \hline
$T$  & Extension parameter
\\ \hline
$N$ & Number of Fourier modes 
\\ \hline
$M$ & Number of equispaced nodes
\\ \hline
$\eta$ & The ratio $M/N$ of nodes to modes
 \\ \hline
 $ \phi_n(x) $ &  The exponential $ \frac{1}{\sqrt{2T}} \E^{\I \frac{n \pi}{T}x}$
 \\ \hline
$\cG^{(T)}_N$ & The space $\spn \{ \phi_n : |n| \leq N \}$
\\ \hline
$S_M$ & The sampling operator $f \mapsto \frac{1}{\sqrt{M}} (f(m/M))^{M}_{m=-M}$
 \\ \hline
 $F^{(T)}_{N,M}(f)$ & The FE approximation \R{FE_LS_fn}
  \\ \hline
 $\mathbf{a} = (a_n)^{N}_{n=-N}$ & Coefficients of the FE approximation 
  \\ \hline
 $A^{(T)}$ & Matrix of the FE approximation
  \\ \hline
   $L^{(T)}_{N,M}$ & Solution of the least-squares \R{LNM_def}
  \\ \hline
$R^{(T)}_{N}$ & The mapping from coefficients to Fourier series     
 \\ \hline
$\epsilon$ & SVD truncation parameter
\\ \hline
$\kappa^{(T)}_{N,M}$ & Condition number
\\ \hline
$\lambda^{(T)}_{N,M}$ & Numerical defect constant
\\ \hline
$E_N(f)$ & Approximation error \R{EN_def}
\\ \hline
$\kappa^*$ & Maximum allowed condition number
\\ \hline
$\Theta^{(T)}(M;\kappa^*)$ & The modes to nodes ratio \R{Theta_def}
\\ \hline
$\nu^{(T)}(\kappa^*)$ & The approximate linear scaling of $\Theta^{(T)}(M;\kappa^*)$ with $M$
\\ \hline
$\tau(\kappa^*)$ & The approximate linear scaling of $\nu^{(T)}(\kappa^*)$ with $T$
\\ \hline
$r$ & The resolution constant  
\\ \hline
\end{longtable}

\end{document}